\documentclass[10pt,a4paper]{article}
\pdfoutput=1 %this serves for arxiv

\usepackage{lineno,hyperref}
\modulolinenumbers[5]

\usepackage{amssymb}
\usepackage{amsmath}
\usepackage{algorithm}     %serve per scrivere lo pseudo codice
\usepackage{algpseudocode} %serve per scrivere lo pseudo codice
\usepackage{subfig}
\usepackage{csquotes}
\usepackage{tikz}
\usetikzlibrary{arrows}
\usepackage{mathtools}
\usepackage{empheq}
\usepackage{bbm}
\usepackage{amsthm}
\usepackage{pgfplots}
\usepackage{float}
\usepackage{cases}
\usepackage{tikz}
\usetikzlibrary{patterns}
\usepackage{xfrac}
\usepackage{setspace}
\usepackage{gensymb}
\usepackage[title]{appendix}
\usepgfplotslibrary{external}
\usepackage{array,booktabs}
\usepackage{geometry}
\newcommand{\email}[1]{\hspace*{\stretch{1}}\emph{\texttt{#1}}}

\makeatletter
\def\blfootnote{\xdef\@thefnmark{$\star$}\@footnotetext}
\makeatother
\newenvironment{Authors}%
  {\begin{center}\begin{bfseries}}%
  {\end{bfseries}\end{center}}
\newenvironment{Addresses}%
  {\begin{flushleft}\begin{itshape}}%
  {\end{itshape}\end{flushleft}}

 \usepackage{fancyhdr}  
  
\fancypagestyle{plain}{
	\fancyhead{}
	\fancyhead[C]{\hfill submitted to International  Journal for Numerical Methods in Engineering, September 2021}
}
\newtheorem{theorem2}{Theorem}[section]

\newtheorem{remark}[theorem2]{Remark}

  \newcommand{\vertiii}[1]{{\left\vert\kern-0.25ex\left\vert\kern-0.25ex\left\vert #1 
    \right\vert\kern-0.25ex\right\vert\kern-0.25ex\right\vert}}

\def\doubleunderline#1{\underline{\underline{#1}}}

\begin{document}

%\begin{frontmatter}
\thispagestyle{plain}

\title{A projection-based model reduction method for nonlinear mechanics with internal variables: application to thermo-hydro-mechanical systems.}
 \date{}
 
 \maketitle
\vspace{-50pt} 
 
\begin{Authors}
Angelo Iollo$^{1}$,
Giulia Sambataro$^{1}$,
Tommaso Taddei$^{1}$,
\end{Authors}

\begin{Addresses}
$^1$
IMB, UMR 5251, Univ. Bordeaux;  33400, Talence, France.
Inria Bordeaux Sud-Ouest, Team MEMPHIS;  33400, Talence, France, \email{angelo.iollo@inria.fr,giulia.sambataro@inria.fr,tommaso.taddei@inria.fr} \\
\end{Addresses}

\begin{abstract}
	We propose a projection-based monolithic model order reduction (MOR) procedure for a class of problems in nonlinear mechanics with internal variables.
	The work is is motivated by applications to thermo-hydro-mechanical (THM) systems for radioactive waste disposal. THM equations model the behaviour of temperature, pore water pressure and solid displacement in the neighborhood of geological repositories, which contain radioactive waste and are responsible for a significant thermal flux  towards the Earth's surface. We develop an adaptive sampling strategy based on the POD-Greedy method, and we develop an element-wise empirical quadrature hyper-reduction  procedure to reduce assembling costs.
	We present numerical results for a   two-dimensional THM system to illustrate and validate the proposed methodology.
\end{abstract}

% REQUIRED
\emph{Keywords:} 
parameterized partial differential equations;  model order reduction; nonlinear mechanics.

%https://onlinelibrary.wiley.com/doi/abs/10.1002/nme.2406?casa_token=c_y-JQhl8_MAAAAA:A_MmIKKVlrNr8MTxPjR8_NKHMLBamOuqvuzC47ypfWASK90LbDc_klDYcL6cvpSzLLsWjSWL3dXXdw
%

\section{Introduction}
\label{sec:introduction}

\subsection{Model reduction for a class of models in nonlinear mechanics}
The disposal and storage of   high-level radioactive waste materials in geological means requires a careful assessment of the long-term effects on neighboring areas.
The system behaviour is well-described by  time-dependent large-scales coupled systems of partial differential equations (PDEs), which take into account the thermal, hydraulic and mechanical response of the geological medium.  Numerical simulation of these systems is 
challenging due to several difficulties:
first, finite element (FE) models of the problem are highly-nonlinear, time-dependent  and high-dimensional; second, due to the uncertainty in model parameters, we need to solve the model for many different system configurations (\emph{many-query} problem).
In this work, we shall devise a model-order reduction (MOR) strategy to speed up parametric studies for radio-active waste disposal applications.

In this contribution we study a general class of nonlinear problems in structural mechanics with internal variables. 
We consider  the spatial variable $x$ in the Lipschitz domain $\Omega \subset \mathbb{R}^d$ with $d=2,3$, and the  time variable $t$ in  the time internal  $(0,T_{\rm{f}})$, where $T_{\rm{f}}$ is the final time. 
We further define the vector of parameters $\mu$ in the compact parameter region   $\mathcal{P} \subset \mathbb{R}^{P}$.
We introduce the state (or primary) variables $\underline{U}$ and  internal (or dependent) variables  $\underline{W}$; 
we denote by $\mathcal{X}$ and $\mathcal{W}$ suitable Hilbert spaces in $\Omega$ for $\underline{U}$ and $\underline{W}$, and we define the space of continuous functions from $(0,T_{\rm{f}})$ to  $\mathcal{X}$ and $\mathcal{W}$, $C(0,T_{\rm{f}}; \mathcal{X})$ and $C(0,T_{\rm{f}}; \mathcal{W})$.
Then, we introduce the parameterised problem of interest: given $\mu\in \mathcal{P}$, find $(\underline{U}_{\mu},\underline{W}_{\mu})\in C(0,T_{\rm{f}}; \mathcal{X}) \times C(0,T_{\rm{f}}; \mathcal{W})$ such that
\begin{equation}
\label{eq:cont_pb}
\left\{
\begin{array}{ll}
\displaystyle{\mathcal{G}_{\mu}(\underline{U}_{\mu},\partial_{\rm{t}}\underline{U}_{\mu}, \underline{W}_{\mu})=0} &
{\rm in} \, \Omega \times (0,T_{\rm{f}}) \\[3mm]
\displaystyle{
\dot{\underline{W}}_{\mu}=\mathcal{F}_{{\mu}}(\underline{U}_{\mu},\underline{W}_{\mu})},
&
{\rm in} \, \Omega \times (0,T_{\rm{f}}) \\
\end{array}
\right.
\end{equation}
with suitable initial and boundary conditions. Here, 
$\mathcal{G}_{\mu}$ is a nonlinear second-order in space, first-order in time differential operator that is associated with the equilibrium equations, while 
$\mathcal{F}_{{\mu}}$ is a set of ordinary  differential equations (ODEs) that is associated with the constitutive laws.

Our  methodology is motivated by the  application to thermo-hydro-mechanical (THM) systems of the form \eqref{eq:cont_pb}, which are widely used to model the system's response for radio-active waste disposal applications.
Radioactive material is placed in an array of  horizontal boreholes (dubbed \emph{alveoli}) deep underground: due to the large temperature of  the alveoli,  a thermal flux is generated; the thermal flux then drives the mechanical and hydraulic response of the medium over the course of several years.
We refer to section \ref{sec:THM_problem} for a detailed discussion of the considered THM model and boundary conditions.

\subsection{Objective of the work and relationship to previous works}
We propose a projection-based monolithic model order reduction (MOR )
(\cite{hesthaven2016certified,quarteroni2015reduced,rozza2007reduced})
technique for problems of the form \eqref{eq:cont_pb}, with particular emphasis on THM systems.
The approach is characterised by an offline/online splitting to reduce the marginal cost, and relies on Galerkin projection to devise a reduced-order model (ROM) for the solution coefficients. We rely on hyper-reduction to speed up the assembly of the ROM during the online stage, and  we rely on adaptive sampling to reduce the offline training costs.

The contribution of the present work is twofold.
First, we propose an element-wise empirical quadrature (EQ) procedure for problems with internal  variables; second, we extend the POD-Greedy algorithm  and we propose  an error indicator that is effective to drive the offline greedy search and is  inexpensive to evaluate.

EQ procedures also dubbed mesh sampling and weighting have been first proposed in
\cite{farhat2015structure,yano2019discontinuous,yano2019lp} and further developed in several other works including \cite{riffaud2021dgdd}: the key feature of EQ is to recast the problem of hyper-reduction as a sparse representation problem and then resort to state-of-the-art  techniques in machine learning and signal processing to estimate the solution to the resulting optimisation problem.
Here, we rely on the approach employed in \cite{taddei2021discretize}, which
combines the methods in 
\cite{farhat2015structure} and \cite{yano2019discontinuous} and 
relies on non-negative least-squares to estimate the solution to the sparse representation problem.
As discussed in section \ref{sec:method}, the presence of internal variables requires several changes to the EQ approach in \cite{taddei2021discretize}.   
We emphasise that several other hyper-reduction techniques have been proposed in the literature including the empirical interpolation method 
(EIM, \cite{barrault2004empirical}) and its discrete variant \cite{chaturantabut2010nonlinear}, the approach in 
 \cite{ryckelynck2009hyper}, and Gappy-POD
\cite{carlberg2013gnat,willcox2006unsteady}: a thorough  comparison of state-of-the-art hyper-reduction techniques is beyond the scope of this work.

The POD-Greedy algorithm was introduced in 
\cite{haasdonk2008reduced} and analysed in \cite{haasdonk2013convergence}: the approach combines proper orthogonal decomposition (POD, \cite{berkooz1993proper,bergmann2009enablers,volkwein2011model}) to compress temporal trajectories with a greedy search driven by an error indicator to explore the parameter domain. In this work, similarly to \cite{fick2018stabilized}, we rely on a time-averaged error indicator to drive the greedy search; furthermore, we test two different compression strategies to update the POD basis at each greedy iteration.

We further observe that the development of online-efficient adaptive ROMs for problems of the form \eqref{eq:cont_pb} is extremely limited in the literature. Relevant examples include the works in 
\cite{ryckelynck2009hyper,miled2013priori,leuschner2017reduced}, which, however, do not consider adaptive sampling.
As regards the application of MOR to THM systems, we recall the recent contribution by Larion \emph{et al.} \cite{larion2020building}: note, however, that the work in  \cite{larion2020building} deals with  a  linearised THM model without internal variables.

\subsection{Outline}
The outline of this paper is the following. In section \ref{sec:formulation} we briefly present the mathematical model  and the numerical discretisation for the general class of nonlinear problems in structural mechanics defined in \eqref{eq:cont_pb}. 
In section \ref{sec:method}, we present the MOR technique: to simplify the presentation, we first discuss the solution reproduction problem and then we extend the approach to the parametric case.
Section \ref{sec:THM_problem} contains details of the THM mathematical model considered in the numerical section.
In section  \ref{sec:numerics},  we present extensive numerical investigations for a two-dimensional THM system. 
In section \ref{sec:conclusions}, we draw some conclusions and we outline a number of subjects of ongoing research.

\section{Formulation}
\label{sec:formulation}
\subsection{Notation}
%Introduco notazioni per dominio, FE spaces, norme, operatore Riesz's, restrizione soluzione agli elementi
In this section, we omit dependence on the parameter.
Given $\Omega\subset \mathbb{R}^d$, we define the triangulation $\{  \texttt{D}_k \}_{k=1}^{N_{\rm e}}$, where $N_{\rm e}$ denotes the total number of elements, the nodes $\{  x_j^{\rm hf}  \}_{j=1}^{N_{\rm hf}}$ and the connectivity matrix $\texttt{T} \in \mathbb{N}^{N_{\rm e}, n_{\rm lp}}$ such that 
$\texttt{T}_{k,i} \in \{1, \ldots,N_{\rm hf} \}$ is the index of the $i$-th node of the $k$-th element of the mesh and $n_{\rm lp}$ is the number of degrees of freedom in each element.  Then, we introduce the  continuous Lagrangian FE basis $\{ \varphi_i  \}_{i=1}^{N_{\rm hf}}$ 
associated with the triangulation $\{  \texttt{D}_k \}_{k=1}^{N_{\rm e}}$, 
such that $\varphi_i(   x_j^{\rm hf}  ) = \delta_{i,j}$,  and we introduce the FE space for the state variables:
\begin{equation}
\label{eq:FE_space}
\mathcal{X}_{\rm hf} := {\rm span}\left\{
\varphi_i \underline{e}_j : 
\;\; i=1,\ldots,N_{\rm hf},
\;\;
j=1,\ldots,D_{\rm eq}
\right\},
\end{equation}
where $\underline{e}_1,\ldots,\underline{e}_{D_{\rm eq}}$ are the elements of the  canonical basis and $D_{\rm eq}$ is the number of state variables. We denote by $\|  \cdot \| = \sqrt{(\cdot,\cdot)}$ the norm of $\mathcal{X}_{\rm hf}$; furthermore, given $\underline{u} \in \mathcal{X}_{\rm hf}$, we denote by $\underline{\mathbf{u}} \in \mathbb{R}^{N_{\rm hf}, D_{\rm eq}}$ the corresponding vector (or matrix) of coefficients such that 
$( \underline{\mathbf{u}}    )_{j,\ell} = \left(   \underline{u}(  x_j^{\rm hf}    )    \right)_{\ell}   $ for $j=1,\ldots,N_{\rm hf}$ and $\ell=1,\ldots,D_{\rm eq}$.

In view of the MOR formulation, we introduce the elemental restriction operators $\mathbf{E}_k : \mathbb{R}^{N_{\rm hf}} \to \mathbb{R}^{n_{\rm lp}}$ such that
\begin{subequations}
\label{eq:restriction_operators}
\begin{equation}
\left( \mathbf{E}_k  \underline{\mathbf{u}}  \right)_{i,\ell}
\, = \,
\left(   \underline{u} (   x_{\texttt{T}_{k,i}}^{\rm hf}    )  \right)_{\ell},
\quad
i=1,\ldots,n_{\rm lp}, \;
\ell=1,\ldots, D_{\rm eq},  \;
k=1,\ldots, N_{\rm e}.
\end{equation}
Furthermore, we introduce the quadrature points 
$\{  x_{q,k}^{\rm hf,q}  \}_{q,k} \subset \Omega$, such that
$ x_{q,k}^{\rm hf,q}$ is the $q$-th quadrature point of the $k$-th  element of the mesh, with $q=1,\ldots,n_{\rm q}$, and the operators
$\mathbf{E}_k^{\rm qd} : \mathbb{R}^{N_{\rm hf}} \to \mathbb{R}^{n_{\rm q}}$  and 
$\mathbf{E}_k^{\rm qd,\nabla} : \mathbb{R}^{N_{\rm hf}} \to \mathbb{R}^{n_{\rm q}, d}$ such that 
\begin{equation}
\left( \mathbf{E}_k^{\rm qd}  \underline{\mathbf{u}}  \right)_{q,\ell}
\, = \,
\left(   \underline{u} (   x_{q,k}^{\rm hf,q}    )  \right)_{\ell},
\quad
\left( \mathbf{E}_k^{\rm qd,\nabla}  \underline{\mathbf{u}}  \right)_{q,\ell,j}
\, = \,
\left(  \frac{\partial }{\partial x_j} \underline{u} (   x_{q,k}^{\rm hf,q}    )  \right)_{\ell}
\end{equation}
where
$q=1,\ldots,n_{\rm q}$,
$\ell=1,\ldots, D_{\rm eq}$,
$k=1,\ldots, N_{\rm e}$ and $j=1,\ldots,d$.
To shorten notation, in the following, we further define
$ \mathbf{E}_k^{\rm qd,\star}: \mathbb{R}^{N_{\rm hf}} \to \mathbb{R}^{n_{\rm q}, d+1}$ such that
\begin{equation}
\left( \mathbf{E}_k^{\rm qd,\star}  \underline{\mathbf{u}}  \right)_{q,\ell,1}
=
\left( \mathbf{E}_k^{\rm qd}  \underline{\mathbf{u}}  \right)_{q,\ell},
\quad
\left( \mathbf{E}_k^{\rm qd,\star}  \underline{\mathbf{u}}  \right)_{q,\ell,2\ldots,d+1}
=
\left( \mathbf{E}_k^{\rm qd,\nabla}  \underline{\mathbf{u}}  \right)_{q,\ell,\cdot}.
\end{equation}
\end{subequations}

\begin{remark}
\label{remark:THM_spoiler}
For the THM problem considered in this work, the state $\underline{U}$ contains the displacement $\underline{u}$, the water pressure $p$ and the temperature $T$ ($D_{\rm eq} = 2+d$); as discussed in the report 
\cite{edfreport2009}, to avoid instabilities, it is important to use polynomials of degree $\kappa$ for displacement and $\kappa-1$ for pressure and temperature: as a result, we should introduce separate restriction operators and separate FE spaces for the different components of the state.
In the main body of the paper we choose to not explicitly address this issue to simplify notation: we remark that the extension to 
$\kappa$-$\kappa-1$ discretisations is computationally tedious but methodologically straightforward.
\end{remark}

\subsection{Finite element  discretisation of \eqref{eq:cont_pb}}
We  introduce the time grid $0=t^{(0)} < t^{(1)} <\ldots < t^{(J_{\rm max})} = T_{\rm f}$ such that $t^{(j)} = j \Delta t$; 
we denote by $\{  \underline{U}_{\rm hf}^{(j)} \}_{j=1}^{J_{\rm max}} \subset \mathcal{X}_{\rm hf}$ the FE approximation of the state variables at each time step; on the other hand, we denote by $\underline{\mathbf{W}}_{\rm hf}^{(j)} \in \mathbb{R}^{n_{\rm q}, N_{\rm e}, D_{\rm int}}$ the tensor associated with the evaluation of the internal variables at time $t^{(j)}$ in the quadrature nodes:
$$
\left( \underline{\mathbf{W}}_{\rm hf}^{(j)}  \right)_{q,k,\ell}
\, = \,
\left(
\underline{W}_{\rm hf}^{(j)} (x_{q,k}^{\rm hf,q})
\right)_{\ell},
\quad
q=1,\ldots,n_{\rm q}, \;
k=1,\ldots,N_{\rm e}.
$$
We further denote by $\texttt{I}_{\rm dir} \subset \{ 1,\ldots,N_{\rm  hf} \cdot D_{\rm eq} \}$ the indices associated with Dirichlet boundary conditions (if any), and we denote by 
$\underline{\mathbf{g}}_{\rm dir}^{(j)} \in \mathbb{R}^{| \texttt{I}_{\rm dir}    |}$ the vector that contains the value of the solution at each Dirichlet node at time $t^{(j)}$.  
%In the remainder of this section, we exploit the duality between combine algebraic notation and

We state the FE discretisation of \eqref{eq:cont_pb} as follows:
for $j=1,2,\ldots$, find $( \underline{U}_{\rm hf}^{(j)}, \underline{W}_{\rm hf}^{(j)}    )$ such that
\begin{equation}
\label{eq:abstract_FE_discretization}
\left\{
\begin{array}{ll}
\displaystyle{
\mathcal{R}^{\rm hf} \left(
\underline{\mathbf{U}}_{\rm hf}^{(j)}, \; \underline{\mathbf{U}}_{\rm hf}^{(j-1)}, \; 
\underline{\mathbf{W}}_{\rm hf}^{(j)}, \; \underline{\mathbf{W}}_{\rm hf}^{(j-1)}, \;  \underline{\mathbf{V}}
\right)
= 0},   \;\;
{\forall \, \underline{V} \in \mathcal{X}_{\rm hf,0}; } \\[3mm]
 \underline{\mathbf{U}}_{\rm hf}^{(j)} (\texttt{I}_{\rm dir}) = \underline{\mathbf{g}}_{\rm dir}^{(j)} ;
\\[3mm]
 \displaystyle{
 \left( \underline{\mathbf{W}}_{\rm hf}^{(j)}  \right)_{q,k,\ell}
 =
\mathcal{F}_{\ell}^{\rm hf} \left(
\left(  \mathbf{E}_k^{\rm qd,\star} \underline{\mathbf{U}}_{\rm hf}^{(j)} \right)_{q,\cdot},  \; \;
 \left(  \mathbf{E}_k^{\rm qd,\star} \underline{\mathbf{U}}_{\rm hf}^{(j-1)} \right)_{q,\cdot},  \; \;
  \left( \underline{\mathbf{W}}_{\rm hf}^{(j-1)}  \right)_{q,k,\cdot}
\right)
   },
\\ [3mm]
\hfill
q=1,\ldots,n_{\rm q}, k = 1,\ldots,N_{\rm e},  \ell=1,\ldots,D_{\rm int}.
\\
 \end{array}
\right.
\end{equation}
where 
$\mathcal{X}_{\rm hf,0}: = \{ \underline{V} \in \mathcal{X}_{\rm hf} \, : \, 
\underline{\mathbf{V}} (\texttt{I}_{\rm dir}) = 0    \} $.
Note that $\mathcal{R}^{\rm hf}$ and $ \underline{\mathcal{F}}^{\rm hf}$ are the discrete counterparts of the operators 
 $\mathcal{G}$ and $ \underline{\mathcal{F}}$ in \eqref{eq:cont_pb}. Note also that the constitutive laws are stated in the quadrature points of the mesh and the internal fields should be computed in the quadrature points of the mesh. 
  
At each time step, 
following \cite{edfreport2009},
we solve \eqref{eq:abstract_FE_discretization} for $\underline{U}_{\rm hf}^{(j)}$ using a Newton method with line search; the method  requires the computation of the Jacobian and the solution to a coupled linear system of size $N_{\rm hf} \cdot D_{\rm eq}$.
Since the underlying problem is second-order in space and first-order in time, the residual $\mathcal{R}^{\rm hf} $ can be written as the sum of local contributions:
\begin{equation}
\label{eq:residual_decomposition}
\begin{array}{l}
\displaystyle{
\mathcal{R}^{\rm hf} \left(
\underline{\mathbf{U}}^{(j)}, \; \underline{\mathbf{U}}^{(j-1)}, \; 
\underline{\mathbf{W}}^{(j)}, \; \underline{\mathbf{W}}^{(j-1)}, \;  \underline{\mathbf{V}}
\right)
\, = \,
}
\\[3mm]
\displaystyle{
\sum_{k=1}^{N_{\rm e}} \; 
r_k^{\rm hf} \left(
\mathbf{E}_{k} \underline{\mathbf{U}}^{(j)}, \; 
\mathbf{E}_{k} \underline{\mathbf{U}}^{(j-1)}, \;
 \left( \underline{\mathbf{W}}^{(j)} \right)_{\cdot,k,\cdot},\;
  \left( \underline{\mathbf{W}}^{(j-1)} \right)_{\cdot,k,\cdot} , \; 
  \mathbf{E}_{k} \underline{\mathbf{V}}^{(j)}
\right)
}
\\
\end{array}
\end{equation}
 As explained in section \ref{sec:method},   decomposition   \eqref{eq:residual_decomposition} provides the foundation of our hyper-reduction procedure.

\section{Methodology}
\label{sec:method}

We propose a time-marching  Galerkin ROM based on linear approximations. More precisely, we consider approximations of the form
\begin{equation}
\label{eq:ROM_ansatz}
\widehat{\underline{U}}_{{\mu}}^{(j)}= \underline{Z} \,  \widehat{\boldsymbol{\alpha}}_{{\mu}}^{(j)}=\sum_{n=1}^{N} \, \left( \widehat{\boldsymbol{\alpha}}_{{\mu}}^{(j)} \right)_n   \underline{\zeta}_n,
\quad
j=1,\ldots,J_{\rm max},
\end{equation}
where $\{  \widehat{\boldsymbol{\alpha}}_{{\mu}}^{(j)}  \}_{j=1}^{J_{\rm max}} \subset \mathbb{R}^N$ are referred to as solution coefficients and are computed by solving a suitable ROM, while $\underline{Z}: \mathbb{R}^N \to \mathcal{X}_{\rm hf}$ is the reduced-order basis (ROB) and
$\mathcal{Z} : = {\rm span} \{ \underline{\zeta}_n \}_{n=1}^N$ is the reduced space.  In presence of non-homogeneous Dirichlet conditions, it is convenient to consider affine approximations of the form 
$\widehat{\underline{U}}_{{\mu}}^{(j)}=   \underline{H} \mathbf{g}^{(j)}  +    \underline{Z} \,  \widehat{\boldsymbol{\alpha}}_{{\mu}}^{(j)}$, where $\underline{H} $ is a suitable lifting operator (see, e.g., \cite{taddei2021discretize}) and $\mathcal{Z} \subset \mathcal{X}_{\rm hf,0}$: since in this work, we consider homogeneous Dirichlet conditions, we do not address the treatment of non-homogeneous  conditions.

The Galerkin ROM is obtained by projecting \eqref{eq:abstract_FE_discretization} onto the reduced space $\mathcal{Z}$: this leads to a nonlinear system of $N$ equations at each time step. To reduce assembly costs, it is important to avoid integration over the whole integration domain. Towards this end, we define the indices associated with  the ``sampled elements'' $\texttt{I}_{\rm eq} \subset \{1,\ldots,N_{\rm e}\}$ and we  define the EQ residual:
\begin{subequations}
\label{eq:galerkin_rom}
\begin{equation}
\label{eq:eq_residuals}
\begin{array}{l}
\displaystyle{
\mathcal{R}_{\mu}^{\rm{eq}} \left(
\underline{\mathbf{U}}^{(j)}, \; \underline{\mathbf{U}}^{(j-1)}, \; 
\underline{\mathbf{W}}^{(j)}, \; \underline{\mathbf{W}}^{(j-1)}, \;  \underline{\mathbf{V}}
\right)
\, = \,}
\\[3mm]
\displaystyle{
\sum_{k \in \texttt{I}_{\rm eq}  }  \;  \rho_k^{\rm eq} \; 
r_{\mu, k}^{\rm hf} \left(
\mathbf{E}_{k} \underline{\mathbf{U}}^{(j)}, \; 
\mathbf{E}_{k} \underline{\mathbf{U}}^{(j-1)}, \;
 \left( \underline{\mathbf{W}}^{(j)} \right)_{\cdot,k,\cdot},\;
  \left( \underline{\mathbf{W}}^{(j-1)} \right)_{\cdot,k,\cdot} , \; 
  \mathbf{E}_{k} \underline{\mathbf{V}}^{(j)}
\right)
}
\\
\end{array}
\end{equation}
where $\boldsymbol{\rho}^{\rm{eq}}=[\rho_1^{\rm{eq}}, ..., \rho_{N_{\rm e}}^{\rm{eq}}]^{T}$ is a sparse vector of positive weights such that 
${\rho}_k^{\rm{eq}} = 0$ if $k\notin \texttt{I}_{\rm eq}$.
In conclusion, the Galerkin ROM reads as follows:
for $j=1,2,\ldots$, find $( \widehat{\underline{U}}_{\mu}^{(j)},  \widehat{\underline{W}}_{\mu}^{(j)}    )$ such that
\begin{equation}
\label{eq:galerkin_rom_b}
\left\{
\begin{array}{ll}
\displaystyle{
\mathcal{R}_{\mu}^{\rm eq} \left(
\widehat{\underline{\mathbf{U}}}_{\mu}^{(j)}, \; \widehat{\underline{\mathbf{U}}}_{\mu}^{(j-1)}, \; 
\widehat{\underline{\mathbf{W}}}_{\mu}^{(j)}, \; 
\widehat{\underline{\mathbf{W}}}_{\mu}^{(j-1)}, \;  \underline{\mathbf{V}}
\right)
= 0},   \;\;
{\forall \, \underline{V} \in \mathcal{Z} ;  } \\[3mm]
 \displaystyle{
 \left( \widehat{\underline{\mathbf{W}}}_{\mu}^{(j)}  \right)_{q,k,\ell}
 =
\mathcal{F}_{\mu, \ell}^{\rm hf} \left(
\left(  \mathbf{E}_k^{\rm qd,\star} \widehat{\underline{\mathbf{U}}}_{\mu}^{(j)} \right)_{q,\cdot},  \; \;
 \left(  \mathbf{E}_k^{\rm qd,\star}\widehat{\underline{\mathbf{U}}}_{\mu}^{(j-1)} \right)_{q,\cdot},  \; \;
  \left( \widehat{\underline{\mathbf{W}}}_{\mu}^{(j-1)}   \right)_{q,k,\cdot}
\right)
  },
\\ [3mm]
\hfill
q=1,\ldots,n_{\rm q}, k\in \texttt{I}_{\rm eq}, \ell=1,\ldots,D_{\rm int}.
\\
 \end{array}
\right.
\end{equation}
Note that the internal variables need to be computed only in the sampled elements. Furthermore, computation of \eqref{eq:galerkin_rom_b} only requires the  storage of the ROB in the sampled elements,
$\{    \mathbf{E}_k \underline{\boldsymbol{\zeta}}_n : n=1,\ldots,N, k\in \texttt{I}_{\rm eq} \}$: provided that $|\texttt{I}_{\rm eq}   | \ll N_{\rm e}$, this leads to significant savings in terms of online assembly costs and also in terms of online memory costs.
\end{subequations}

In  the remainder of this section, we shall discuss the construction of the ROB $\underline{Z}$ (\emph{data compression}), the empirical quadrature rule 
$\boldsymbol{\rho}^{\rm{eq}}$ (\emph{hyper-reduction}) and also the  error indicator. To simplify the presentation, 
in section \ref{sec:solution_reproduction_pb} we focus on  the solution reproduction problem, while in section \ref{sec:parametric_pb} we discuss the extension to the parametric problem.

\subsection{Solution reproduction problem}
\label{sec:solution_reproduction_pb}

The  solution reproduction problem refers to the task of reproducing the results obtained for a fixed value of the parameter $\mu$.  Algorithm \ref{alg:SRP}
summarises the procedure:
during the offline stage, we compute the hf solution to \eqref{eq:ROM_ansatz} for a given parameter and we store snapshots of the state variables at 
select time steps  $\texttt{I}_{\rm s} \subset \{1,\ldots,J_{\rm max} \}$;
then, we use this piece of information to build a ROM for the state;
then, during the online stage, we query the ROM for the same value of the parameter considered in the offline stage.

The solution reproduction problem is of little practical interest; however, it represents  the first step towards
the implementation of an effective  ROM for the parametric problem. Note  that during the offline stage we store the state variables in a subset of the time steps and we do not store internal variables: this choice is motivated by the fact that for practical problems memory constraints might prevent the storage of all snapshots; in addition, internal variables might not be computed explicitly by available hf codes.

\begin{algorithm}[H]
\caption{Solution reproduction problem: offline/online decomposition}
\begin{algorithmic}[1]
\Statex{\textbf{Offline stage:}}
\State compute 
$\{  \underline{\mathbf{U}}_{\rm hf,\mu}^{(j)}  \}_{j\in \texttt{I}_{\rm s}}$, 
$\texttt{I}_{\rm s} \subset \{1,\ldots,J_{\rm max} \}$;
\medskip

\State 
construct the ROB $\underline{Z}$;
 \Comment{section \ref{sec:POD}}  
 \smallskip

\State 
construct the weights $\boldsymbol{\rho}^{\rm{eq}}$.
 \Comment{section \ref{sec:eq}}  
 \smallskip
 
\Statex\textbf{Online stage:}
\State compute $\{ \widehat{\boldsymbol{\alpha}}_{\mu}^{(j)}  \}_{j=1}^{J_{\rm max}}$  by  solving  the ROM \eqref{eq:galerkin_rom}.
	\end{algorithmic}
\label{alg:SRP}
\end{algorithm}

\subsubsection{Data compression}
\label{sec:POD}
 
We resort to POD based on the method of snapshots  (cf. \cite{sirovich1987turbulence}) to generate the ROB $\underline{Z}$.
Given the snapshots
$\{  \underline{{U}}_{\rm hf,\mu}^{(j)}  \}_{j\in \texttt{I}_{\rm s}}
=
\{  \underline{{U}}^{(k)}  \}_{k=1}^K$, $K = |  \texttt{I}_{\rm s} |$, we define the Gramian matrix $\mathbf{C} \in \mathbb{R}^{K,K}$  such that
$\mathbf{C}_{k,k'}=(\underline{U}^k,\underline{U}^{k'})$; then, we define the  POD eigenpairs
$$
\mathbf{C} \widetilde{\boldsymbol{\zeta}}_n = \lambda_n  \widetilde{\boldsymbol{\zeta}}_n,
\quad
\lambda_1 \geq \lambda_2\geq \ldots \lambda_K \geq 0;
$$
finally, we define the POD modes
$$
\underline{\zeta}_n : =
\sum_{k=1}^K \; \left( \widetilde{\boldsymbol{\zeta}}_n  \right)_k \underline{U}_k,
\quad
n=1,\ldots,N.
$$

The  reduced space size $N$  can be chosen according to the energy criterion:
\begin{equation}
\label{eq:POD_maxN}
N:=\min \biggl\{M:  \sum_{n=1}^{M} \lambda_{n} 
\geq 
(1-tol_{\rm{POD}}^2)  \; \sum_{i=1}^{K} \; \lambda_i \biggr\},
\end{equation}
for some user-defined tolerance $tol_{\rm{POD}}>0$.
Note that the POD modes depend on the choice of the inner product $(\cdot, \cdot)$: we discuss the choice of $(\cdot, \cdot)$ for the THM problem considered in this paper in section \ref{sec:THM_problem}.

\subsubsection{Hyper-reduction}
\label{sec:eq}

We denote  by $\widehat{\mathbf{R}}_{\mu}^{\rm hf}(\cdot)$ and 
$\widehat{\mathbf{R}}_{\mu}^{\rm eq}(\cdot)$ the algebraic reduced residuals associated with the hf and empirical quadrature rules, such that
$$
\left\{
\begin{array}{l}
\displaystyle{
 \left( \widehat{\mathbf{R}}_{\mu}^{\rm hf} \left( \boldsymbol{\alpha}; \,    \boldsymbol{\beta},  \underline{\mathbf{W}}'  \right) \right)_n
: =
\mathcal{R}_{\mu}^{\rm hf}
\left(
\underline{\mathbf{Z}}  \, \boldsymbol{\alpha},  \; 
\underline{\mathbf{Z}}  \, \boldsymbol{\beta}, \; 
\underline{\mathbf{W}}_{\mu}^{\star}, \; 
\underline{\mathbf{W}}', \;  \underline{\boldsymbol{\zeta}}_n
\right), 
\quad
n=1,\ldots,N, 
}
\\[4mm]
\displaystyle{
 \left( \widehat{\mathbf{R}}_{\mu}^{\rm eq} \left( \boldsymbol{\alpha};    \boldsymbol{\beta},  \underline{\mathbf{W}}' \right) \right)_n
: =
\mathcal{R}_{\mu}^{\rm eq}
\left(
\underline{\mathbf{Z}}  \, \boldsymbol{\alpha},  \; 
\underline{\mathbf{Z}}  \, \boldsymbol{\beta}, \; 
\underline{\mathbf{W}}_{\mu}^{\star} , \; \underline{\mathbf{W}}', \;  \underline{\boldsymbol{\zeta}}_n
\right), 
\quad
n=1,\ldots,N, 
}
\\
\end{array}
\right.
$$
where $\boldsymbol{\alpha}, \boldsymbol{\beta} \in \mathbb{R}^N$,
$\underline{\mathbf{W}}' \in \mathbb{R}^{n_{\rm q}, N_{\rm e}, D_{\rm int}}$, 
 and 
$\underline{\mathbf{W}}_{\mu}^{\star} = \underline{\mathbf{W}}_{\mu}^{\star}\left( \boldsymbol{\alpha},  \boldsymbol{\beta}; \underline{\mathbf{W}}'\right)$ is obtained by substituting in \eqref{eq:abstract_FE_discretization}$_3$.
We further introduce the Jacobians
$\mathbf{J}_{\mu}^{\rm hf}(\cdot),\mathbf{J}_{\mu}^{\rm eq}(\cdot)$ such that
$$
\left( \mathbf{J}_{\mu}^{\rm hf}(  \boldsymbol{\alpha};    \boldsymbol{\beta},  \underline{\mathbf{W}}'    ) \right)_{n,n'}
: =
\frac{\partial}{\partial \alpha_{n'} } \left( \widehat{\mathbf{R}}_{\mu}^{\rm hf} \left( \boldsymbol{\alpha};    \boldsymbol{\beta}, \underline{\mathbf{W}}'  \right) \right)_n,
\quad
\left( \mathbf{J}_{\mu}^{\rm eq}(  \boldsymbol{\alpha};    \boldsymbol{\beta},   \underline{\mathbf{W}}'    ) \right)_{n,n'}
: =
\frac{\partial}{\partial \alpha_{n'} } \left( \widehat{\mathbf{R}}_{\mu}^{\rm eq} \left( \boldsymbol{\alpha};    \boldsymbol{\beta},  \underline{\mathbf{W}}'  \right) \right)_n,
$$
for $n,n'=1,\ldots,N$. We observe that the computation of the Jacobian involves the derivatives with respect to the constitutive laws in $\underline{\mathcal{F}}^{\rm hf}$; we further observe that the residuals
$\widehat{\mathbf{R}}_{\mu}^{\rm hf}(\cdot)$ and 
$\widehat{\mathbf{R}}_{\mu}^{\rm eq}(\cdot)$ satisfy
\begin{equation}
\label{eq:boring_identity}
\widehat{\mathbf{R}}_{\mu}^{\rm hf} \left( \boldsymbol{\alpha};    \boldsymbol{\beta}, \underline{\mathbf{W}}'  \right)  =
\mathbf{G} \left( \boldsymbol{\alpha};    \boldsymbol{\beta},  \underline{\mathbf{W}}'  \right)  \; \boldsymbol{\rho}^{\rm hf}, 
\quad
\widehat{\mathbf{R}}_{\mu}^{\rm eq} \left( \boldsymbol{\alpha};    \boldsymbol{\beta},  \underline{\mathbf{W}}'  \right)  =
\mathbf{G} \left( \boldsymbol{\alpha};  \boldsymbol{\beta}, \underline{\mathbf{W}}'  \right)  \; \boldsymbol{\rho}^{\rm eq}, 
\end{equation}
where $\mathbf{G} \in \mathbb{R}^{N,N_{\rm e}}$ can be explicitly derived using the same approach  as in \cite{taddei2021discretize} and
$\boldsymbol{\rho}^{\rm hf} = [1,\ldots,1]^T$.

As in  \cite{yano2019lp}, we reformulate the problem of finding the sparse weights $\boldsymbol{\rho}^{\rm eq} \in \mathbb{R}^{N_{\rm e}}$ as the problem of finding a vector 
$\boldsymbol{\rho}^{\rm eq}$ such that:
\begin{enumerate}
\item 
the number of nonzero entries in $\boldsymbol{\rho}^{\rm eq}$, 
which we denote by $\|  \boldsymbol{\rho}^{\rm eq} \|_0$,  is as small as possible;
\item 
the entries of $\boldsymbol{\rho}^{\rm eq}$  are non-negative;
\item 
(\emph{constant-function constraint})
the constant function is integrated accurately:
$\Big|\displaystyle{\sum_{k=1}^{N_{\rm e}}}\rho_k^{\rm{eq}} |\texttt{D}_k|  -|\Omega|\Big| \ll 1$;
\item
(\emph{manifold accuracy constraint})
the empirical and hf residuals are close at operating conditions:
\begin{equation}
\label{eq:manifold_accuracy}
\big\|
\left( \mathbf{J}_{\mu}^{\rm hf}(  \boldsymbol{\alpha}_{\rm train}^{(j)},    \boldsymbol{\alpha}_{\rm train}^{(j)}; \underline{\mathbf{W}}_{\rm train}^{(j-1)}    ) \right)^{-1}
\left(
\widehat{\mathbf{R}}_{\mu}^{\rm hf} \left(  \boldsymbol{\alpha}_{\rm train}^{(j)},    \boldsymbol{\alpha}_{\rm train}^{(j)}; \underline{\mathbf{W}}_{\rm train}^{(j-1)}    \right)
\; -  \;
\widehat{\mathbf{R}}_{\mu}^{\rm eq} \left(  \boldsymbol{\alpha}_{\rm train}^{(j)},    \boldsymbol{\alpha}_{\rm train}^{(j)}; \underline{\mathbf{W}}_{\rm train}^{(j-1)}    \right)
\right)
\big\|_2 \ll 1,
\end{equation}
for $j\in \texttt{I}_{\rm s}$ and for suitable choices of 
$\{  \boldsymbol{\alpha}_{\rm train}^{(j)} \}_j$ and 
$\{  \underline{\mathbf{W}}_{\rm train}^{(j)}  \}_j$ that are discussed at the end of the section.
\end{enumerate}

Exploiting \eqref{eq:boring_identity}, we can restate the previous requirements as a sparse representation problem:
\begin{equation}
\label{eq:hyper_sparse_pb1}
\rm{find} \, \boldsymbol{\rho}^{\rm{eq}} \in \arg \displaystyle{\min_{\boldsymbol{\rho}\in \mathbb{R}^{N_{\rm e}}}} 
\| \boldsymbol{\rho} \|_0 \, \rm{s.t.} \, 
\begin{cases}
 \boldsymbol{\rho}\geq \mathbf{0} \\
\|   \mathbf{C} \boldsymbol{\rho}-\mathbf{b} \|_{*}\leq \delta,
\end{cases}
\end{equation}
for a suitable choices of the matrix $\mathbf{C}$, the  vector $\mathbf{b}$, 
the norm $\|   \cdot  \|_{*}$, 
and the tolerance $\delta$. Since \eqref{eq:hyper_sparse_pb1} is NP-hard, we find an approximate solution
to \eqref{eq:hyper_sparse_pb1}
 by solving the non-negative least-squares problem:
\begin{equation}
\label{eq:hyper_sparse_pb}
\displaystyle{\min_{\boldsymbol{\rho}\in \mathbb{R}^{N_{\rm e}}}} \| 
\mathbf{C}  \boldsymbol{\rho}-\mathbf{b}||_2 \,  \rm{s.t.} \, \boldsymbol{\rho}\geq \mathbf{0}.
\end{equation}
In this work, we rely on the Matlab function \texttt{lsqnonneg} that implements the Greedy algorithm proposed in \cite{lawson1974solving} and takes as input the matrix $\mathbf{C}$, the  vector $\mathbf{b}$, and a tolerance $tol_{\rm eq}$:
$$
\boldsymbol{\rho}^{\rm{eq}}=\mathtt{lsqnonneg}(\mathbf{C} , \mathbf{b}, tol_{\rm eq}).
$$
The same algorithm to find the sparse weights 
$\boldsymbol{\rho}^{\rm{eq}}$ given the matrices 
$\mathbf{C}, \mathbf{b}$ has been first considered in \cite{farhat2015structure}: for large-scale problems, a parallelised extension  of the algorithm was introduced and successfully applied to hyper-reduction in \cite{chapman2017accelerated}.

We remark that in order to compute the entries of  $\mathbf{C} , \mathbf{b}$ associated with \eqref{eq:manifold_accuracy}  we should prescribe the triplets
$\left\{ \left(  \boldsymbol{\alpha}_{\rm train}^{(j)},
\boldsymbol{\alpha}_{\rm train}^{(j-1)},
\underline{\mathbf{W}}_{\rm train}^{(j-1)} 
\right) \right \}_{j\in \texttt{I}_{\rm s}}$;
note, in particular, that the internal variables cannot be directly extracted from   hf computations. We here choose to consider
$ \boldsymbol{\alpha}_{\rm train}^{(j)} =  
\widehat{\boldsymbol{\alpha}}_{\rm hf, \mu}^{(j)}$ and
$\underline{\mathbf{W}}_{\rm train}^{(j)} =  
\widehat{\underline{\mathbf{W}}}_{\rm hf, \mu}^{(j)}$
where 
$\{ \widehat{\boldsymbol{\alpha}}_{\rm hf, \mu}^{(j)}, \widehat{\underline{\mathbf{W}}}_{\rm hf, \mu}^{(j)} \}_j$ denote the solution to \eqref{eq:galerkin_rom} for 
$\boldsymbol{\rho}^{\rm eq} =\boldsymbol{\rho}^{\rm hf}$.  
Note that this choice  requires the solution to a ROM with hf quadrature.

\subsection{Parametric problem}
\label{sec:parametric_pb}

In order to extend our methodology to parametric problems, we should address two challenges. First, we should propose an adaptive strategy to explore the parameter domain $\mathcal{P}$ based on an inexpensive error indicator; second, we should devise a compression strategy to combine information from different parameters.

Our point of departure is the POD-Greedy algorithm proposed in  \cite{haasdonk2008reduced}. Algorithm \ref{alg:HPODGreedy} summarises the procedure: 
the procedure takes as input a discretisation of $\mathcal{P}$,
$\Xi_{\rm train}$, a tolerance $tol_{\rm loop}$ for the outer greedy loop, a tolerance $tol_{\rm pod}$ for the data compression step, and the maximum number of greedy iterations $N_{\rm count, max}$ --- we here prescribe the termination condition based on the error indicator; we refer to the pMOR literature for other termination conditions.

We observe that the algorithm depends on several building blocks: the FE solver
$$
\left[ \{  \underline{\mathbf{U}}_{\rm hf,\mu}^{(j)}  \}_{j\in \texttt{I}_{\rm s}} \right]
= \texttt{FE-solve} ( \mu )
$$
takes  as input the vector of parameters and returns the snapshot set associated with the sampling times
$\texttt{I}_{\rm s} \subset \{1, \ldots, J_{\rm max} \}$;
the data compression routine
$$
\left[  \underline{Z}', \; \boldsymbol{\lambda}' \right]
= \texttt{data-compression} \left(  \underline{Z}, \; \boldsymbol{\lambda} ,
\{  \underline{\mathbf{U}}_{\rm hf,\mu^{\star}}^{(j)}  \}_{j\in \texttt{I}_{\rm s}} ,
(\cdot,\cdot), tol_{\rm pod} \right)
$$
takes as input the current ROB and the POD eigenvalues 
$\boldsymbol{\lambda} = [\lambda_1,\ldots,\lambda_N]^T$, and returns the  updated  ROB $\underline{Z}'$ and the updated  eigenvalues $\boldsymbol{\lambda}'$;
finally, we observe that construction of the ROM comprises both the construction of the Galerkin ROM and of the error indicator. 
In the remainder of this section, we discuss each element of the procedure.

\begin{algorithm}[htb]
\setstretch{1.3}
\caption{POD-Greedy algorithm}
\begin{algorithmic}[1]
\Require 
$\Xi_{\rm train} = \{ \mu^{(k)}  \}_{k=1}^{n_{\rm train}}$,
 $tol_{\rm loop}$, $tol_{\rm pod}$, $N_{\rm count, max}$.
\smallskip

\State
$\mathcal{Z} =\emptyset$,
$\boldsymbol{\lambda} =\emptyset$,
 $\mu^{\star}=\mu^{(1)} $.
 \smallskip
 
\For {$n_{\rm{count}} = 1,\ldots,  N_{\rm count, max}$}
 
\State 
$\left[ \{  \underline{\mathbf{U}}_{\rm hf,\mu^{\star}}^{(j)}  \}_{j\in \texttt{I}_{\rm s}} \right]
= \texttt{FE-solve} ( \mu^{\star} )$;
\smallskip

\State 
$\left[  \underline{Z}, \; \boldsymbol{\lambda} \right]
= \texttt{data-compression} (\underline{Z}, \; \boldsymbol{\lambda} ,
\{  \underline{\mathbf{U}}_{\rm hf,\mu^{\star}}^{(j)}  \}_{j\in \texttt{I}_{\rm s}} ,
(\cdot,\cdot), tol_{\rm pod} )$;
 \Comment{section  \ref{sec:parametric_data_compression}.}
\smallskip

\State 
Construct the ROM with error indicator.
 \Comment{section  \ref{sec:parametric_ROM_construction}.}
\smallskip
 
\For {$j=1:n_{\rm{train}}$}
\State 
Solve the ROM \eqref{eq:galerkin_rom} for $\mu= \mu^{(k)}$ and compute $\Delta_{\mu}$.		
\EndFor
\smallskip

\State 
$\mu^{\star}=\arg \max_{ {\mu} \in \Xi_{\rm{train}}} \Delta_{\mu}$ \Comment{Greedy search}
\smallskip

\If {$\Delta_{\mu^{\star}}< tol_{\rm loop}$},
\Comment{Termination condition}

\State 
\texttt{break}, \EndIf.

\EndFor

\Return 
ROB $\underline{Z}$ and ROM: $\mu\in \mathcal{P} \mapsto 
\{ \widehat{\boldsymbol{\alpha}_{\mu}}^{(j)}  \}_{j=1}^{J_{\rm max}}$.

\end{algorithmic}
\label{alg:HPODGreedy}
\end{algorithm}

\subsubsection{Data compression}
\label{sec:parametric_data_compression}

We consider two different data compression strategies: a hierarchical POD (H-POD) and a hierarchical approximate POD (HAPOD). Both techniques have been considered in several previous works:
we refer to \cite[section 3.5]{haasdonk2017reduced} for H-POD and to 
\cite{himpe2018hierarchical} for HAPOD; HAPOD is also related to incremental singular value decomposition in linear algebra \cite{brand2003fast}. Here, we review the two approaches for completeness. We denote by $\Pi_{\mathcal{Z}}: \mathcal{X}_{\rm hf} \to \mathcal{Z}$ the orthogonal projection operator on $\mathcal{Z} \subset \mathcal{X}_{\rm hf} $; furthermore, we introduce notation
$$
\left[  \underline{Z}, \; \boldsymbol{\lambda} \right]
= \texttt{POD} \left(  
\{  \underline{U}^{(k)}  \}_{k=1}^K ,
(\cdot,\cdot), tol_{\rm pod}  \right)
$$
to refer to the application of POD to the snapshot set 
$\{  \underline{U}^{(k)}  \}_{k=1}^K$, with inner product $(\cdot, \cdot)$, and tolerance $tol_{\rm pod}$
(cf. \eqref{eq:POD_maxN}), with 
$\underline{Z} = [ \underline{\zeta}_1,\ldots,\underline{\zeta}_N]$,
$\| \underline{\zeta}_n \| = 1$,
$\boldsymbol{\lambda} =
[\lambda_1,\ldots,    \lambda_N]^T$, and
$\lambda_1\geq \lambda_2 \ldots \geq \lambda_N$.

Given  $\underline{Z}$ and the snapshots
$\{  \underline{U}_{\rm hf, \mu^{\star}}^{(j)}   \}_{j}$
H-POD considers the update:
\begin{subequations}
\label{eq:HPOD}
\begin{equation}
\label{eq:HPOD_a}
\underline{Z}'=
\left[
\underline{Z}, 
\underline{Z}^{\rm new}
\right],
\quad
 \underline{Z}^{\rm new}
= \texttt{POD} \left(  
\{  
\Pi_{\mathcal{Z}^{\perp}}  \underline{U}_{\rm hf, \mu^{\star}}^{(j)}
 \}_{j} ,
(\cdot,\cdot), tol_{\rm pod}  \right).
\end{equation}
Note that the approach does not require to input the POD eigenvalues $\boldsymbol{\lambda}$ from the previous iterations. We observe that the approach leads to a sequence of nested spaces --- that is, the updated ROB contains the ROB of the previous iteration --- and  it returns an orthonormal basis of the reduced space. In our experience, the choice of the tolerance $tol_{\rm pod}$ is extremely challenging: since 
\eqref{eq:POD_maxN} depends on the relative energy content of the snapshot set, 
the update \eqref{eq:HPOD_a} with fixed tolerance $tol_{\rm pod}$ 
might lead  to an excessively large (resp., small) number of modes when 
$\max_j \| \underline{U}_{\rm hf, \mu^{\star}}^{(j)}  - 
\Pi_{\mathcal{Z}}  \underline{U}_{\rm hf, \mu^{\star}}^{(j)}        \|$ is small (resp., large). For this reason, we propose to choose the number of new modes $N^{\rm new}$ using the criterion:
\begin{equation}
N^{\rm new} \, := \,\min 
\left\{
M \,: \,
\max_{j \in \texttt{I}_{\rm s}} \;
\frac{ 
\|  \Pi_{( \mathcal{Z}\oplus \mathcal{Z}_{M}^{\rm new} )^{\perp}}  \underline{U}_{\rm hf, \mu^{\star}}^{(j)} \| }      
{ \|    \underline{U}_{\rm hf, \mu^{\star}}^{(j)} \| }
\leq tol_{\rm pod},
\;\;
\mathcal{Z}_{M}^{\rm new} = {\rm span} \{ 
\underline{\zeta}_m^{\rm new}
\}_{m=1}^M
\right\}.
\label{eq:HAPOD_N}
\end{equation}
Note that this choice enforces that the in-sample relative projection error is below a certain threshold for all snapshots computed during the greedy iterations.
\end{subequations}

HAPOD considers the update
\begin{equation}
\label{eq:HAPOD}
[\underline{Z}', \boldsymbol{\lambda}' ]    
= \texttt{POD} \left(  
\{    \underline{U}_{\rm hf, \mu^{\star}}^{(j)}  \}_{j}
\cup \{ \lambda_n \underline{\zeta}_n \}_{n=1}^N ,
(\cdot,\cdot), tol_{\rm pod}  \right).
\end{equation}
Note that the approach \eqref{eq:HAPOD} does  not in general  lead to hierarchical (nested) spaces.
As discussed in \cite[section 3.3]{himpe2018hierarchical}, 
which refers to \eqref{eq:HAPOD} as to \emph{distributed HAPOD}, 
it is possible to relate the performance of the reduced space obtained using 
HAPOD  to the performance of the
POD space associated with the snapshot set
$\{  
 \underline{U}_{\rm hf, \mu^{\star,n}}^{(j)} \, : \,
 n=1,\ldots,N_{\rm count,max} ,
 \;
 j\in \texttt{I}_{\rm s} \}$: we refer to  the above-mentioned paper for a thorough discussion.

\subsubsection{Time-averaged error indicator}
\label{sec:parametric_error_indicator}
%\overrightarrow{ \underline{U} }, \overrightarrow{ \underline{W} },
We define the trajectories
$\mathbb{U} = \{   \underline{U}^{(j)} \}_{j=1}^{J_{\rm max}} $ and
$\mathbb{W} = \{   \underline{W}^{(j)} \}_{j=1}^{J_{\rm max}} $; given the pair $\left(  \mathbb{U} , \mathbb{W}    \right)$, we define the time-average residual:
\begin{equation}
\label{eq:time_avg_residual}
\mathcal{R}_{\rm avg, \mu}^{\rm hf} \left(
 \mathbb{U} , \mathbb{W},
\underline{V}
\right) \, := \,
\sum_{j=1}^{J_{\rm max}} \;
(t^{(j)} - t^{(j-1)}  )
\; 
\mathcal{R}_{\mu}^{\rm hf} \left(
\underline{U}^{(j)} ,  \underline{U}^{(j-1)} , 
\underline{W}^{(j)} ,  \underline{W}^{(j-1)} ,
\underline{V}
\right),\quad
\forall \; \underline{V} \in \mathcal{X}_{\rm hf,0},
\end{equation}
and the error indicator
\begin{equation}
\label{eq:hf_dual_residual}
\Delta_{\mu}^{\rm hf} \left(  
 \mathbb{U} , \mathbb{W} 
 \right)
 \; = \;
 \sup_{  \underline{V} \in \mathcal{X}_{\rm hf,0}      }
 \;\;
 \frac{ \mathcal{R}_{\rm avg, \mu}^{\rm hf} \left(
 \mathbb{U} , \mathbb{W},  \underline{V} \right)}{
\| \underline{V}  \|}.
\end{equation}
The indicator \eqref{eq:hf_dual_residual} is expensive to evaluate since it relies on hf quadrature and it requires the computation of the supremum over all elements of $\mathcal{X}_{\rm hf,0}$: following  \cite{taddei2019offline}, we consider the  hyper-reduced error indicator
\begin{equation}
\label{eq:error_indicator}
\Delta_{\mu}  \left(   \mathbb{U} , \mathbb{W}  \right)
 \; = \;
 \sup_{  \underline{V} \in \mathcal{Y}    }
 \;\;
 \frac{ \mathcal{R}_{\rm avg, \mu}^{\rm eq,r} \left(
\mathbb{U} , \mathbb{W}, \underline{V} \right)}{
\| \underline{V}  \|},
\end{equation}
where $\mathcal{Y} \subset \mathcal{X}_{\rm hf,0}$ is an $M$-dimensional empirical test space, while 
$\mathcal{R}_{\rm avg, \mu}^{\rm eq,r}$ is defined by replacing 
 $\mathcal{R}_{\mu}^{\rm hf}$ in \eqref{eq:time_avg_residual} with a suitable sparse weighted residual of the form \eqref{eq:eq_residuals}, defined over the elements $\texttt{I}_{\rm eq,r}\subset \{1,\ldots,N_{\rm e}\}$.
 
 Given the ROM solution $ ( 
\widehat{\mathbb{U}}_{\mu} ,  \widehat{\mathbb{W}}_{\mu}
)$, the test space $\mathcal{Y}$ should guarantee that
\begin{equation}
\label{eq:error_desiderata_calY}
 \sup_{  \underline{V} \in \mathcal{Y}    }
 \;\;
 \frac{ \mathcal{R}_{\rm avg, \mu}^{\rm hf} \left(
\widehat{\mathbb{U}}_{\mu} ,  \widehat{\mathbb{W}}_{\mu},
\underline{V} \right)}{
\| \underline{V}  \|} \approx
 \sup_{  \underline{V} \in \mathcal{X}_{\rm hf,0}    }
 \;\;
 \frac{ \mathcal{R}_{\rm avg, \mu}^{\rm hf} \left(
\widehat{\mathbb{U}}_{\mu} ,  \widehat{\mathbb{W}}_{\mu},
\underline{V} \right)}{
\| \underline{V}  \|},
\quad
\forall \mu \in \mathcal{P},
\end{equation}
which implies that $ \mathcal{Y}  $ should be an approximation of the space of Riesz elements 
$\mathcal{M}_{\rm test} : = \{  \widehat{\underline{\psi}}_{\mu} \, : \, \mu\in  \mathcal{P} \}$ with 
\begin{equation}
\label{eq:riesz_representers}
\left(
\widehat{\underline{\psi}}_{\mu}, \underline{V}
\right) \, =\,
\mathcal{R}_{\rm avg, \mu}^{\rm hf} \left(
\widehat{\mathbb{U}}_{\mu} ,  \widehat{\mathbb{W}}_{\mu},
\underline{V} \right), \quad
\forall \; \underline{V} \in \mathcal{X}_{\rm hf,0}.
\end{equation}
 On  the other hand, the empirical quadrature rule should ensure that
 \begin{equation}
 \label{eq:accuracy_constraints_residual}
 \mathcal{R}_{\rm avg, \mu}^{\rm eq,r} \left(
\widehat{\mathbb{U}}_{\mu} ,  \widehat{\mathbb{W}}_{\mu}, 
\underline{\psi}_m \right) \, \approx  \,
 \mathcal{R}_{\rm avg, \mu}^{\rm hf} \left(
\widehat{\mathbb{U}}_{\mu} ,  \widehat{\mathbb{W}}_{\mu}, 
\underline{\psi}_m \right),
\quad
\forall \, \mu \in \mathcal{P}, \;
m=1,\ldots,M,
 \end{equation}
 where $\underline{\psi}_1,\ldots,\underline{\psi}_M$ is an orthonormal basis of $\mathcal{Y}$. 
 
In our implementation, we compute the error indicator during the time iterations --- as opposed to after having computed the whole solution trajectory.  Algorithm \ref{alg:error_indicator_computation} provides the complete online solution and residual indicator computations.
We find that computation of $\Delta_{\mu}$ requires to compute the internal variables 
 $\widehat{\mathbb{W}}_{\mu}$ in the elements
 $\texttt{I}_{\rm eq}\cup \texttt{I}_{\rm eq,r}$  at each time iteration (cf. \eqref{eq:galerkin_rom_b}), and it requires to store the trial ROB $\underline{Z}$ in  $\{\texttt{D}_k : k\in    \texttt{I}_{\rm eq}\cup \texttt{I}_{\rm eq,r}\}$ and the test basis $\underline{Y} = [ \underline{\psi}_1,\ldots,\underline{\psi}_M   ]$ in $\{\texttt{D}_k : k\in \texttt{I}_{\rm eq,r}\}$.

\begin{algorithm}[htb]
\setstretch{1.3}
\caption{Online solution and residual computations}
\begin{algorithmic}[1]

\State
Initial state and internal variables; set 
$ \widehat{\mathbf{R}}_{\mu}^{\rm avg} = \mathbf{0}$.
\smallskip

\For {$j=1,\ldots, J_{\rm max}$}
\State
Compute $\widehat{\boldsymbol{\alpha}}_{\mu}^{(j)}$ by solving \eqref{eq:galerkin_rom_b}.
 \smallskip
 
\State 
Compute
$\left( \widehat{\underline{\mathbf{W}}}_{\mu}^{(j)} \right)_{\cdot, k, \cdot}$ for all $k\in \texttt{I}_{\rm eq,r}$ using \eqref{eq:galerkin_rom_b}$_2$.
\smallskip

\State 
Assemble
$  \widehat{\mathbf{R}}_{\mu}^{(j)} \in \mathbb{R}^M$ 
such that
$\left(
\widehat{\mathbf{R}}_{\mu}^{(j)} 
\right)_m = 
 \mathcal{R}_{\mu}^{\rm eq,r} \left(
{\widehat{\underline{U}}}_{\mu}^{(j)} , \; 
{\widehat{\underline{U}}}_{\mu}^{(j-1)} , \; 
{\widehat{\underline{W}}}_{\mu}^{(j)} , \;
{\widehat{\underline{W}}}_{\mu}^{(j-1)}, \;
\underline{\psi}_m \right)$ for $m=1,\ldots,M$.
\smallskip

\State 
Update
$  \widehat{\mathbf{R}}_{\mu}^{\rm avg} =   
\widehat{\mathbf{R}}_{\mu}^{\rm avg}  + (t^{(j)} - t^{(j-1)}) \widehat{\mathbf{R}}_{\mu}^{(j)}$.
\smallskip
\EndFor

\Return 
$\{  \widehat{\boldsymbol{\alpha}}_{\mu}^{(j)}   \}_j$ and
$\Delta_{\mu} = \|  \widehat{\mathbf{R}}_{\mu}^{\rm avg}  \|_2$

\end{algorithmic}
\label{alg:error_indicator_computation}
\end{algorithm}
  
  Several authors  (e.g., \cite{haasdonk2008reduced}) have considered the time-discrete
$L^2(0,T_{\rm f}; \mathcal{X}_{\rm hf,0}')$ 
residual indicator
\begin{equation}
\label{eq:L2_residual}
\Delta_{\mu}^{\rm hf,2} \left(  
 \mathbb{U} , \mathbb{W} 
 \right)
 \; = \;
\sqrt{
\sum_{j=1}^{J_{\rm max}} \;
(t^{(j)} - t^{(j-1)}  )
\; 
 \left(
 \sup_{  \underline{V} \in \mathcal{X}_{\rm hf,0}      }
 \;\;
 \frac{ 
\mathcal{R}_{\mu}^{\rm hf} \left(
\underline{U}^{(j)} ,  \underline{U}^{(j-1)} , 
\underline{W}^{(j)} ,  \underline{W}^{(j-1)} ,
\underline{V}
\right) 
 }{
\| \underline{V}  \|}
\right)^2
}. 
\end{equation}
We observe that  we could apply the same ideas considered in this section to devise an hyper-reduced counterpart of the residual indicator \eqref{eq:L2_residual}. However, we find that the test space $\mathcal{Y}$ and the empirical quadrature rule should be accurate for all parameters and for all time steps: as a result, the resulting test space $\mathcal{Y}$ might be significantly higher dimensional  and the quadrature rule might be significantly less sparse, for the desired accuracy. For this reason, in this work, we investigate the effectivity of the time-averaged error indicator \eqref{eq:error_indicator}.

\subsubsection{ROM construction}
\label{sec:parametric_ROM_construction}

In order to devise an actionable ROM, we should discuss 
(i) the choice of the EQ rule $\boldsymbol{\rho}^{\rm eq}$, 
(ii) the choice of the test space $\mathcal{Y}$ and   of the EQ rule $\boldsymbol{\rho}^{\rm eq,r}$ in \eqref{eq:error_indicator}. In view of the presentation of the computational procedure, we 
define the ROM solution with hf quadrature
$( \widehat{\mathbb{U}}_{\mu} ^{\rm hf} , 
 \widehat{\mathbb{W}}_{\mu} ^{\rm hf})$;
we denote by $\mathbf{C}_{\mu} \in \mathbb{R}^{K\cdot N, N_{\rm e}}$ the EQ matrix 
associated with the manifold accuracy constraints in \eqref{eq:manifold_accuracy}   for $\mu\in \mathcal{P}$ (cf. section \ref{sec:eq});  
we further define the vector 
$\mathbf{c} = [|\texttt{D}_1|, \ldots, | \texttt{D}_{N_{\rm e}}  |]^T$ associated with the constant function accuracy constraint.
Given the  test reduced basis
$\underline{\psi}_1,\ldots,\underline{\psi}_M$, 
we define $\mathbf{G}_{\mu}^{\rm r} \in \mathbb{R}^{M,N_{\rm e}}$ such that
\begin{equation}
\label{eq:Gmat}
\left(  \mathbf{G}_{\mu}^{\rm r}   \boldsymbol{\rho}_{\rm hf} 
\right)_m
\, = \,
 \mathcal{R}_{\rm avg, \mu}^{\rm eq,r} \left(
\widehat{\mathbb{U}}_{\mu}^{\rm hf} ,  \widehat{\mathbb{W}}_{\mu}^{\rm hf}, \; 
\underline{\psi}_m \right),
\quad
\forall \; \mu\in  \mathcal{P}, \;
m=1,\ldots,M.
\end{equation}
We further define the unassembled average residual 
$\mathbf{R}_{\mu}^{\rm avg, un} \in \mathbb{R}^{n_{\rm lp}, N_{\rm e}, D_{\rm eq}}$:  we observe that 
$\mathbf{R}_{\mu}^{\rm avg, un}$ might be employed to build the FE residual and ultimately compute the Riesz representers $\widehat{\underline{\psi}}_{\mu}$ in \eqref{eq:riesz_representers}, and also, given $\mathcal{Y}$, to compute 
$\mathbf{G}_{\mu}^{\rm r}$.

We focus on the construction of the ROM at the $n_{\rm c}$-th iteration of the POD Greedy algorithm. We define
$\Xi^{\star} = 
\{\tilde{\mu}^{(j)} \}_{j=1}^{n_{\rm rom}} = 
\{ \mu^{\star, (i)}  \}_{i=1}^{n_{\rm c}} \cup 
 \{ \tilde{\mu}^{(j)}  \}_{j=1}^{n_{\rm train,eq}} 
$,  where
$\mu^{\star, (1)}, \ldots,\mu^{\star, (n_{\rm c})}$ are the parameters sampled by the greedy algorithm and
$\tilde{\mu}^{(1)}, \ldots, \tilde{\mu}^{(n_{\rm train,eq})}$
are independent identically distributed samples from  the uniform distribution over $\mathcal{P}$. Algorithm \ref{alg:ROM_construction} summarises the computational procedure as implemented in our code. The test space $\mathcal{Y}$ is built using POD as in \cite{taddei2019offline}, while the EQ weights $\underline{\rho}_{\rm{eq},r}$ are obtained using the non-negative least-squares method.

\begin{algorithm}[htb]
\setstretch{1.3}
\caption{Construction of the ROM}
\begin{algorithmic}[1]

\For {$\mu \in \Xi^{\star} $}
\State
Solve the ROM with hf quadrature and compute 
$\mathbf{C}_{\mu}$ and $\mathbf{R}_{\mu}^{\rm avg, un}$.
 \EndFor
\smallskip 
 
\State 
Assemble
$
\mathbf{C} = \left[
\begin{array}{l}
\mathbf{C}_{\tilde{\mu}^{(1)} } \\
\vdots \\
\mathbf{C}_{\tilde{\mu}^{(n_{\rm rom})} } \\
\mathbf{c}^T \\
\end{array}
\right] \in \mathbb{R}^{K\cdot N\cdot n_{\rm rom}, N_{\rm e}}$
and set
$\boldsymbol{\rho}^{\rm{eq}}=\mathtt{lsqnonneg}(\mathbf{C} , \mathbf{C} \boldsymbol{\rho}^{\rm{hf}}, \, tol_{\rm eq})$.
\smallskip 
 
\State 
Compute the Riesz representers $\{  \widehat{\underline{\psi}}_{\mu}  \}_{\mu \in \Xi_{\star}}$ using \eqref{eq:riesz_representers}.
\smallskip 
 
\State 
Define the empirical test space 
$\mathcal{Y} = {\rm span} \{ \underline{\psi}_m \}_{m=1}^M$ as
$[ \{ \underline{\psi}_m \}_{m=1}^M ] 
= \texttt{POD} \left(  
\{  \widehat{\underline{\psi}}_{\mu}  \}_{\mu \in \Xi_{\star}} ,
(\cdot,\cdot), tol_{\rm pod,res}  \right).$
\smallskip
 
 \State
Assemble
$
\mathbf{G} = \left[
\begin{array}{l}
\mathbf{G}_{\tilde{\mu}^{(1)} } \\
\vdots \\
\mathbf{G}_{\tilde{\mu}^{(n_{\rm rom})} } \\
\mathbf{c}^T \\
\end{array}
\right]  \in \mathbb{R}^{M \cdot n_{\rm rom}, N_{\rm e}}$
and set
$\boldsymbol{\rho}^{\rm{eq,r}}=\mathtt{lsqnonneg}(\mathbf{G} , \mathbf{G} \boldsymbol{\rho}^{\rm{hf}}, \, tol_{\rm eq,r})$.
 
\end{algorithmic}
\label{alg:ROM_construction}
\end{algorithm}

\section{The THM model}
\label{sec:THM_problem}
 
 In this section we illustrate the non-dimensional mathematical formulation and the numerical discretisation of the THM system considered in this work. We assume  that the solid undergoes small displacements and that soil is fully-saturated in water.
We resort to a Lagrangian formulation for the solid, and to an Eulerian formulation for the fluid.
%(since we deal with classic mechanics, without dissipative or friction terms)

\subsection{Preliminary definitions}
\label{sec:THM_definitions}

We first introduce the state variables and the internal variables. 
The state variables represent solid displacement, water pressure and temperature and are reported in Table \ref{table:primary}; the internal variables $\underline{W}=[\rho_{\rm{w}}, \varphi, h_{\rm{w}}, Q, \underline{M}_{\rm{w}}^{\rm{T}}, m_{\rm{w}}]^{\rm{T}}$
represent dependent physical quantities  and are illustrated in Table \ref{table:internal}, together with the corresponding SI units.

\begin{table}[H]
	\centering
	\begin{tabular}{lll}
		\toprule
		& SI unit      & description              \\
		\midrule
		$\underline{u}$ & $\rm{m}$  & solid displacement \\
		$p_{\rm{w}}$    & $\rm{Pa}$ & water pressure     \\
		$T$             & $\rm{K}$  & temperature       \\
		\hline
	\end{tabular}
	\caption{primary variables}
	\label{table:primary}
\end{table}

\begin{table}[H]
	\centering
	\begin{tabular}{lllll}
		\toprule
		& SI unit                     & label                 \\ 
		\midrule
		$\rho_{\rm{w}}$  & $\rm{kg \cdot m^{-3}}$ & water density \\
		$\varphi$ & $\%$    & Eulerian porosity  \\
		$h_{\rm{w}} $     &   $\rm J \cdot Kg^{-1}$         & mass enthalpy of water  \\ 
		$\mathcal{Q}$        & $\rm Pa$& non-convected heat\\
		$ \underline{M}_{\rm{w}}$    &  $\rm kg \cdot m^{-2} \cdot s^{-1}$      &     mass flux \\
		$m_{\rm{w}} $         &$\rm kg \cdot m^{-3}$& mass input       \\
		\hline
	\end{tabular}
	\caption{dependent variables}
	\label{table:internal}
\end{table}

We denote the Cauchy stress tensor  by $ \doubleunderline{\sigma}  [\rm Pa]$, and 
we define the volumetric deformation $ \epsilon_{\rm V}=tr(\doubleunderline{\epsilon})$ where $\doubleunderline{\epsilon}$ is the strain tensor: $\doubleunderline{\epsilon}=\nabla_{\rm s} \underline{u} = \frac{1}{2}\left( \nabla \underline{u} +  \nabla \underline{u}^{T}\right)$ . 
We also provide 
in Table \ref{table:characteristic_par}
  the characteristic parameters that
we use for the non-dimensionalisation.

\begin{table}[H]
	\centering
	\begin{tabular}{llll}
		\toprule
		& SI unit                                  & value               \\
		\midrule
		$\bar{t}$             & $\rm s$                                             & $3.15 \cdot 10^{7}$               \\
		$\bar{H}$             & $\rm m$               & $77.3$              \\
		$\sigma_0$            & $\rm Pa$                                            & $11.3 \cdot 10^{6}$ \\
		$\rho_0$              & $\rm kg \cdot m^{-3}$                             & $2450$              \\
		$T_{\rm ref}$         & $\rm K$                                            & $297.5$             \\
		$\overline{\Delta T}$ & $\rm K$                                          & $30$  \\              
		\hline
	\end{tabular}
	\caption{characteristic constants}
	\label{table:characteristic_par}
\end{table}

\subsubsection{Geometry configuration}
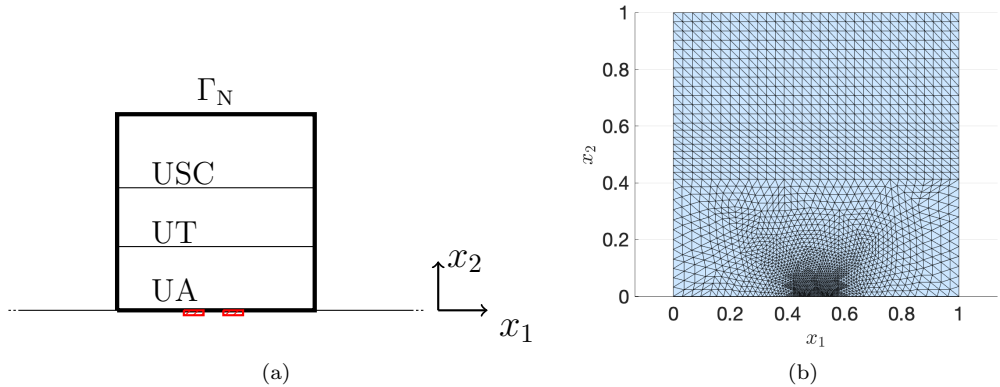
\begin{figure}[H]
	\centering
	\subfloat[]{\label{fig:domain}
		\begin{tikzpicture}[scale=.65]
			\linethickness{0.3 mm}
			\linethickness{0.3 mm}
			
			\draw[ultra thick]  (0,0)--(4,0)--(4,4)--(0,4)--(0,0);
			\draw[-]  (-2,0)--(6,0);
			\draw[-]  (6.05,0)--(6.1,0);
			\draw[-]  (6.15,0)--(6.2,0);
			
			\draw[-]  (-2.10,0)--(-2.05,0);
			\draw[-]  (-2.2,0)--(-2.15,0);
			\draw[-]  (0,1.3)--(4,1.3);
			\draw[-]  (0,2.5)--(4,2.5);

			\draw[thick,red,pattern=north east lines,pattern color=red]    (1.35,0)-- (1.75,0)-- (1.75,-0.1)-- (1.35,-0.1)-- (1.35,0);
			
			\draw[thick,red,pattern=north east lines,pattern color=red]    (2.15,0)-- (2.55,0)-- (2.55,-0.1)-- (2.15,-0.1)-- (2.15,0);
			
			\coordinate [label={right:  {\large {UA}}}] (E) at (0.5, 0.42) ;
			\coordinate [label={right:  {\large {UT}}}] (E) at (0.5, 1.6) ;
			
			\coordinate [label={right:  {\large {USC}}}] (E) at (0.5, 2.8) ;
			\draw[thick, ->]  (6.5,0)--(6.5,1);
			\coordinate [label={below:  {\Large {$x_1$}}}] (E) at (8.1,0) ;
			\draw[thick, ->]  (6.5,0)--(7.5,0);
			\coordinate [label={left:  {\Large {$x_2$}}}] (E) at (7.6,1) ;
			
			\coordinate [label={above:  {\large {$\Gamma_{\rm N}$}}}] (E) at (2, 4) ;
		\end{tikzpicture}
	%\label{fig:domain}
	}
	\subfloat[]{\label{fig:geometry} \includegraphics[width=.38\textwidth]{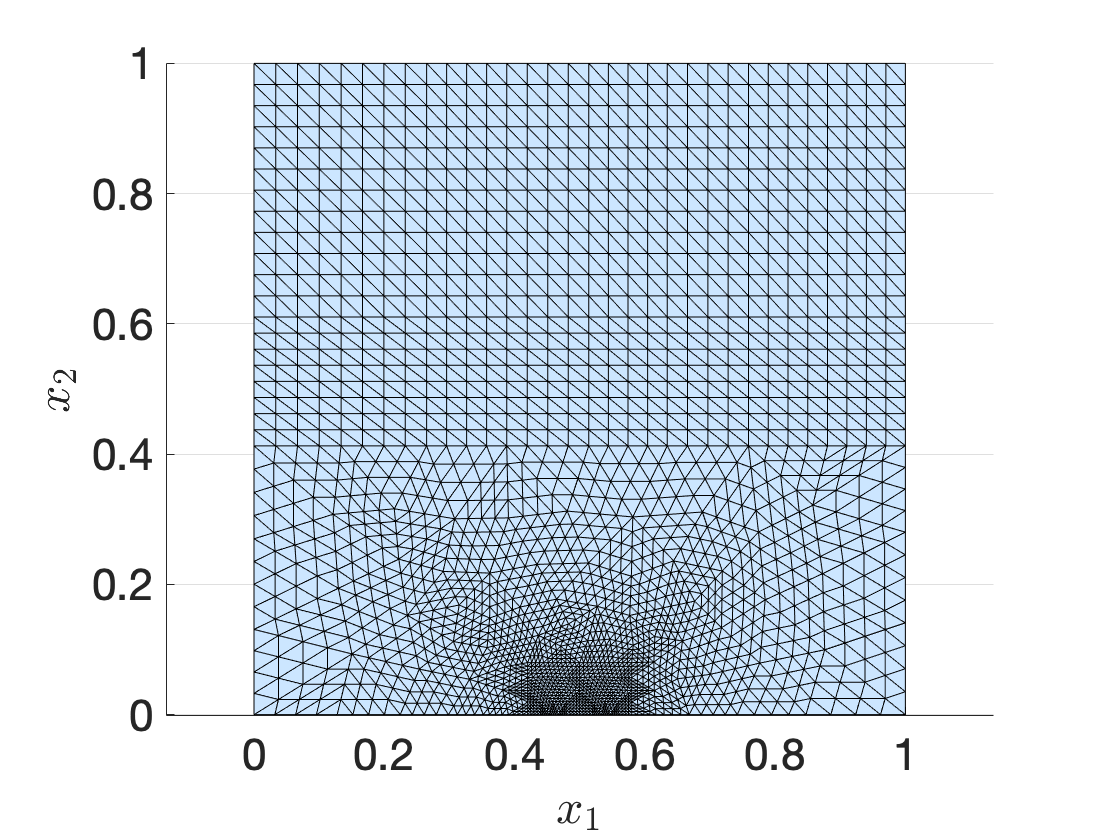}\label{fig:mesh}}
	\caption{geometric configuration: (a) the non-dimensional domain (b) the mesh. The size of each alveoulus is equal to $l_{Q}=3.09 \, [\rm m]$, while the distance between consecutive alveoli is equal to $l=6.18 \, [\rm m].$}
	\label{fig:geometry_mesh}
\end{figure}

The computational domain is shown in Figure~ \ref{fig:geometry_mesh}\subref{fig:domain}. The geological repositories, modelled as boundary conditions, are depicted in red at the bottom of the domain,  in the case of two activated alveoli. In  the vertical $(x_2)$ direction, the domain is split into three layers: a clay layer denoted as UA ("\textit{unit\'e argilleuse}"), a transition layer UT ("\textit{unit\'e de transition}") and a silt-carbonate layer USC ("\textit{unit\'e silto-carbonat\'ee}").   \\
In Figure~ \ref{fig:geometry_mesh}\subref{fig:mesh} the finite element grid is shown. The number of degrees of freedom for the first state component (solid displacement) is $N_{\rm{hf}}^{\rm{u}}=40430$, while for water pressure and temperature is $N_{\rm{hf}}^{\rm{p}}=N_{\rm{hf}}^{\rm{t}}=9045$. 
The grid is refined in the proximity of the alveoli to better capture the relevant features of the solution.
We consider  a $\texttt{p}=3$ FE discretisation for the displacement component, and a $\texttt{p}=2$ FE discretisation for both pressure and temperature.
\subsection{Mathematical problem}
We first state the equilibrium equations -- the superscripts $(\cdot)^{\rm{m}}$, $(\cdot)^{\rm{n}}$,$(\cdot)^{\rm{t}}$ refer to quantities associated with the mechanical, hydraulic and thermal behaviours, respectively. Then, we present the constitutive laws that are considered and finally we present the boundary conditions. To clarify the presentation, we report in Table \ref{table:further} the parameters that enter in the constitutive laws.

We denote by $\underline{F}_{\rm{m}}=-\frac{g}{\gamma} \underline{e}_{2}$ (where $\gamma=\frac{\sigma_0}{\rho_0 \bar{H}}$) the mechanical force with $g$ defined in Table \ref{table:further} and we specify that $\underline{n}$ (resp. $\underline{t}$) is the unitary outward normal (resp. tangential) vector in the domain depicted in Figure \ref{fig:geometry_mesh}\subref{fig:domain}; then we introduce the equilibrium of mechanical forces:
\begin{subequations}
	\begin{equation}
		\left\{
		\begin{array}{ll}
			-\nabla \cdot \doubleunderline{\sigma} \, =\, \rho \underline{F}_{\rm m} & {\rm in} \, \Omega, \\[3mm]
			\doubleunderline{\sigma}\, \underline{n} = \underline{g}_{\rm m,N}
			& {\rm on} \, \Gamma_{\rm N},\\[3mm]
			\underline{u} \cdot \underline{n}  = 0
			& {\rm on} \, \partial\Omega \setminus \Gamma_{\rm N},\\[3mm]
			(\doubleunderline{\sigma}\, \underline{n})\cdot \underline{t}=0  & {\rm on} \, \partial\Omega \setminus \Gamma_{\rm N},
		\end{array}
		\right.
		\label{eq:mechanics}
	\end{equation}
where $\Gamma_{\rm N}$ is depicted in Figure~ \ref{fig:geometry_mesh}\subref{fig:domain}.
	The Neumann datum $\underline{g}_{\rm m,N}$ is given by $\underline{g}_{\rm m,N}=-\underline{e}_2$. The stress tensor is linked to the primary and internal variables by the linear law
	\begin{equation}
		\doubleunderline{\sigma}= 
		2 {\mu} {\nabla}_{\rm s} {\underline{u}}
		+ \left(
		{\lambda}  \, {\nabla} \cdot  {\underline{u}} \,  - \, 
		(2{\mu} + 3 {\lambda}) {\alpha_{\rm s}} {T} 
		- b {p}_{\rm w}  \right)   \mathbbm{1}, 
	\end{equation}
\end{subequations}
where the Lam\'e constants $\mu', \, \lambda$ satisfy
\begin{align*}
	&\mu'=\frac{E}{2(1 + \nu)} ,\\
	& \lambda=\frac{E \nu}{(1+\nu)(1-2 \nu)},
\end{align*} 
and $E$ and $\nu$ are introduced in Table \ref{table:further}.

We state the mass conservation of water as follows
\begin{subequations}
	\begin{equation}
		\left\{
		\begin{array}{ll}
			\partial_{{t}} \, {m}_{\rm w} \, + \, {\nabla} \cdot {\underline{M}}_{\rm w} = 0
			& {\rm in} \, \Omega \\[3mm]
			{\underline{M}}_{\rm w} \cdot  \underline{n} =  0
			& {\rm on} \, \partial \Omega \\ 
		\end{array}
		\right.
		\label{eq:hydraulics}
	\end{equation}
	where the muss flux $\underline{M}_{\rm{w}}$ is given by the Darcy law
	\begin{equation}
		\displaystyle{
			{\underline{M}}_{\rm w} = - {\gamma} \left(   
			{\nabla}  {p}_{\rm w}  - {\rho}_{\rm w} {\underline{F}}_{\rm m} \right)},
	\end{equation}
	and 
	\begin{equation}
		\displaystyle{
			{\gamma} = 	\rho_{\rm{w}}\frac{\kappa_{\rm w} \, \sigma_0 \, \bar{t}}{\rho_0 \mu_{\rm w,0} \, \bar{H}^2 }  \, 
			{\rm exp}\left(
			-\frac{1808.5}{T_{\rm ref} + \overline{\Delta T} \, {T}  }
			\right)}.
	\end{equation}
\end{subequations}

Finally we consider the energy balance:
\begin{subequations}
	\begin{equation}
		\left\{
		\begin{array}{ll}
			{h}_{\rm w} \partial_{{t}} \, {m}_{\rm w} \, + \, 
			\partial_{{t}} \, {\mathcal{Q}} \, + \, 
			{\nabla} \cdot \left(  {h}_{\rm w}  {\underline{M}}_{\rm w} \, + \, {\underline{q}}  \right) - {\underline{M}}_{\rm w} \cdot {\underline{F}}_{\rm m} \,=\, {\Theta}
			& {\rm in} \, \Omega \\[3mm]
			\left(  {h}_{\rm w}  {\underline{M}}_{\rm w} \, + \, {\underline{q}} \right) 
			\cdot  \underline{n} =  {g}_{\rm t, N}
			& {\rm on} \, \partial \Omega \\ 
		\end{array}
		\right.
		\label{eq:energy}
	\end{equation}
	where $\mathcal{Q}$ is the non-convective heat, $\underline{q}$ is the thermal flux and is given by the Fick law 
	\begin{equation}
		\displaystyle{
			{\underline{q}}= - {\Lambda} {\nabla} {T},
		}
	\end{equation}
	with ${\Lambda}=\rm{diag}({\lambda}_{1},{\lambda}_{2})$.
\end{subequations}
If we denote by $\Gamma_{\rm{al}} \subset \partial \Omega$ the region associated with the alveoli, $g_{t,N}$ is equal to 
\begin{equation}
	g_{t,N}=\frac{P_{\rm{t}}n_c \bar{t}}{l_{Q}\bar{H}^2 \sigma_0}\exp\big(-t/\tau\big)  \mathbbm{1}_{\Gamma_{\rm al}}=C_{\rm{al}} \exp\big(-t/\tau\big)\mathbbm{1}_{\Gamma_{\rm{al}}},
	\label{eq:gtN}
	\end{equation}
where $n_c \,  [\%]$ is the density of the radioactive waste stock in each alveolus (equal to $45$ anisters), $P_{\rm{t}}=31.4 \, [\rm{W}]$ is the unitary termic power at the initial time, $l_{Q}=3.09 \, [\rm{m}]$ is the size of each alveolus, $\sigma_0, \bar{H}, \bar{t}$ are introduced in Table \ref{table:characteristic_par} and $\tau=\frac{\bar{t}}{\log(0.112)} \, [\rm{s}]$ is a characteristic decay time.

\begin{subequations}
	\begin{empheq}[left=\empheqlbrace]{align}
		\label{eq:rhow}
		&\displaystyle{
			\frac{d\rho_{\rm w}}{ \rho_{\rm w}} \, = \, 
			\frac{d p_{\rm w}}{K_{\rm w}} \, - \, 
			3 {\alpha}_{\rm w} dT}\\
		&	\displaystyle{
			\frac{d \varphi}{b \, - \, \varphi}  \, = \, 
			d \epsilon_{\rm V} \, - 3 {\alpha}_{\rm s} d{T} + \frac{d {p}_{\rm w}}{{ K}_{\rm s}}},\\
		&\displaystyle{
			d {h}_{\rm w}  = {C}_{\rm w}^{\rm p} \, d {T} \, + \, ( \beta_{h}^{\rm p} - 3 {\alpha}_{\rm w} {T} ) \frac{d {p}_{\rm w}}{ {\rho}_{\rm w}}	}, \\
		& \displaystyle{
			\delta {\mathcal{Q}} = 
			\left( \beta_{\mathcal{Q}}^{\epsilon}  + 3 {\alpha_{\rm s}} { K}_0 \, {T}  \right) \, d \epsilon_{\rm V} 
			\, -
			\left( \beta_{\mathcal{Q}}^{\rm p}  + 
			3 {\alpha}_{\rm w,m}   {T}  \right) \, d {p}_{\rm w}
			\, + \, {C}_{\epsilon}^0 \, d{T}},\\
		&	\displaystyle{{m}_{\rm w} = {\rho}_{\rm w}  (1 + \epsilon_{\rm V}) \,  \varphi  - {\rho}_{\rm w}^0 \varphi^0}
		\label{eq:mw}
	\end{empheq}	
\end{subequations}
Here, we have  $\beta_{h}^{p}=1-3 \alpha_{\rm{w}}T_{\rm ref}$, $\beta_{\mathcal{Q}}^{\epsilon}=3 \alpha_{\rm s} K_0 T_{\rm ref}$, $\beta_{\mathcal{Q}}^{p}=3 \alpha_{\rm w,m} T_{\rm{ref}}$.\\ The parameters in \eqref{eq:rhow}-\eqref{eq:mw} are defined in Table \ref{table:further}.

\begin{table}[H]
	\small
	\begin{tabular}{lllll}
		\toprule
		&SI  unit                	& description           	& reference value			&  formula \\
		\midrule
		$g$ 	&  	$\rm m \cdot s^{-2}$   & gravity acceleration	 & $9.81$           &        \\
		$E$       &      $\rm Pa$  						 &  Young's modulus            &  \begin{tabular}[c]{@{}l@{}}
			$11.4 \cdot 10^9 \: \rm UA$\\ $12.3\cdot 10^9   \: {\rm UT}$\\ $20 \cdot 10^9    \: {\rm USC}$
		\end{tabular}
		& 	\\ 
		$\nu$&              $\%$                                         &  Poisson's ratio     &    $0.3$             &                     \\
		$\mu$ & $\rm Pa$ 							& Lam\'e parameter,		    & 									&$\frac{E}{2(1 + \nu)} $ \\
		$\lambda$ & $\rm Pa$ 						 & Lam\'e parameter		&    						&  		$\frac{E \nu}{(1+\nu)(1-2 \nu)}$								\\
		$b$ 					& $\%$  					&Biot coefficient					& $0.6$					& \\[4.5mm]
		$\alpha_{\rm s}$ & $\rm K^{-1}$ & solid thermal expansion coefficient & $1.28 \cdot 10^{-5}$ & \\
		$\alpha_{\rm 0}$& $\rm K^{-1}$ & expansion coefficient											& $1.28 \cdot 10^{-5}$ & \\
		$\kappa_{\rm w}$ & $\rm m^2$ & intrinsic permeability of  porous medium 			& $10^{-21}$					&   \\
		$\mu_{\rm w}$& $\rm MPa \cdot s$ &   dynamic viscosity											& 								&$\mu_{\rm w}=\mu_{\rm w,0} \exp (\frac{1808.5}{T})$ \\
		$\mu_{\rm w,0}$ & $\rm MPa \cdot s$ & dynamic viscosity coefficient					& $2.1 \cdot 10^{-12} $ & \\
		$K_{\rm s}$ 			& $\rm Pa$ 					& bulk modulus of the solid								&												&	\\
		$K_{\rm{w}}$	& $\rm Pa$ & 				bulk modulus of water									&		$ 2\cdot 10^9$				&		$K_{\rm s}=\frac{E}{3(1-2\nu)} $	\\
		$C_{\rm w}^{\rm p}$& $\rm J \cdot kg^{-1} \cdot K^{-1}$	& 	heat capacity at constant pressure	&	$4180$										&	\\[4.5mm]
		$K_0$& 	$\rm Pa$	& 			drained bulk modulus	&		&$K_0=(1-b)K_{\rm s} $	\\
		$\alpha_{\rm{w}}$& 	$K^{-1}$			& 	thermal expansion coefficient of water				&							&	$\alpha_{\rm{w}}=9.52 \cdot10^{-5} \log(T-273)-2.19 \cdot 10^{-4}$\\
		$\alpha_{\rm w,m}$& 						& 		dilation coefficient	&					&	\\
		$C_{\sigma}^{\rm s}$& $\rm J \;  kg^{-1}  \cdot K$								& 			specific heat at constant stress				&	 
		\begin{tabular}[c]{@{}l@{}}
			$537 \: {\rm UA}$	\\
			$603 \: {\rm UT}$	 \\ 
			$640  \: {\rm USC}$	 
		\end{tabular}
		&	\\
		$\rho^0$& $\rm Kg \cdot m^{-3}$& porous medium initial density&
		\begin{tabular}[c]{@{}l@{}}
			$2450 \: {\rm UA}$\\
			$2450  \: {\rm UT} $\\
			$2500  \: {\rm USC}$
		\end{tabular}
		&\\
		$\rho_{\rm{w}}^0$&$\rm Kg \cdot m^{-3}$					&initial water density &$10^{3}$	& 		\\
		$\varphi^0$& $\%$	& 	initial Eulerian porosity		&		
		\begin{tabular}[c]{@{}l@{}}
			$0.25 \: {\rm UA}$\\
			$0.21  \: {\rm UT}$\\ 
			$0.19  \: {\rm USC}$
		\end{tabular}				&	\\
		$h_{\rm w}^0$ &  $\rm m^2 \cdot s^{-2}$ &  initial water enthalpy&   &$h_{\rm w}^0=\frac{p_{\rm w}^0-p_{\rm atm}}{\rho_{\rm w}^0}$\\
		$\rho_{\rm s}$& $\rm Kg \cdot m^{-3}$										& 			density ratio							&										&$\rho_{\rm s}=\frac{\rho^0 -\rho_{\rm w}^0 \varphi^0}{1-\varphi^0}$	\\
		$C_{\epsilon}^0$&$\rm Pa \cdot K^{-1}$ 		& 	specific heat at constant deformation								&													&$C_{\epsilon}^0=(1-\varphi)\rho_s C_{\sigma}^s +\varphi \rho_w C_w^p-9 T K_{0} \alpha_{\rm s}^2$	\\
		$\Lambda$& 					& 	thermic conductivity tensor							&					&			$\Lambda=\rm{diag}(\lambda_1, \lambda_2)$ 	\\
		$\lambda_1$& $\rm W m^{-1}K^{-1}$& thermic conductivity component& $\begin{array}{ll}
			1.5 & {\rm UA} \\
			1.5 & {\rm UT} \\
			1.3 & {\rm USC} \\
		\end{array}$&\\
		$\lambda_2$&  $\rm W m^{-1}K^{-1}$&thermic conductivity component & $\begin{array}{ll}
			1 & {\rm UA} \\
			1 & {\rm UT} \\
			1.3 & {\rm USC} \\
		\end{array}$&\\
		$\Theta$& $\rm Pa \cdot s^{-1}$& volumetric heat sources& &\\
		\hline
	\end{tabular}
	\caption{parameters of the constitutive laws. Layers UA, UT, USC are depicted in Figure~\ref{fig:geometry_mesh} \protect\subref{fig:domain}.}
	\label{table:further}
\end{table}

\subsubsection{Initial conditions}
To set the initial conditions, we consider the case of deactivated repositories: therefore, we set thermal flux equal to zero and we set a constant temperature $T_0=T_{\rm{ref}}$ in $\Omega$, where the reference temperature is defined in Table \ref{table:characteristic_par}.
 We aim at finding the initial values of the primary variables $\underline{u}$ and $p_{\rm{w}}$ that correspond to the equilibium solutions of a preliminary problem: here, the Neumann boundary condition for the energy equation is zero, that is, $g_{\rm t,N}=0$, and temperature is costant and equal to the reference value $T_{\rm{ref}}$ (in Table \ref{table:characteristic_par}).\\
 We then seek $\underline{u}_0$, $p_{\rm{w},0}$ such that the initial solution vector $\underline{U}_0=[\underline{u}^{\rm{T}}_0, p_{\rm{w}}, T_0]^{\rm{T}}$ satisfies the equilibrium equations \eqref{eq:mechanics}, \eqref{eq:hydraulics} and \eqref{eq:energy} with thermal flux $g_{\rm{t,N}}$ equal to $0$ on the  domain boundary $\partial \Omega$. Towards this end, we first observe that \eqref{eq:rhow} reduces to 
 \begin{equation}
 	\frac{d\rho_{\rm{w}}}{\rho_{\rm{w}}}=\frac{dp_{\rm{w}}}{K_{\rm{w}}}
 \end{equation}
that brings to
 $p_{\rm{w}}=\rho_{-\infty}\exp\left(\frac{1}{K_{\rm{w}}}(p_{\rm{w}}-p_{-\infty})\right)$. If we assume that $\rho_{\rm{w}}=\rho_{-\infty}=\rho_{\rm{w},0}$, we find $p_{\rm{w}}=p_{-\infty}$; furthermore, by susbstituting these assumptions into the hydraulic equilibrium equation we find
 \begin{equation}
 	p_{\rm{w},0}(x, y)=p_{\rm{w, top}}+\rho_{\rm{w},0}g (1-y)
 \end{equation}
where $p_{\rm{w, top}}$ is a datum for water pressure that is defined at the top boundary of the domain $(0,1)\times \{1\}$.
Finally, we search for $\underline{u}_0$ as the solution to the  equilibrium equation of mechanical forces:
\begin{equation}
	\int_{\Omega} \,  2 \mu \, \nabla_{\rm s} \, \underline{u}_0 \, : \,  \nabla_{\rm s} \, \underline{v} \, + \,\lambda (\nabla  \cdot \underline{u}_{0}) (\nabla \cdot \underline{v})- b p_{\rm{w},0}   \, \nabla \cdot \underline{v} \, - \, {\rho}^0  \underline{F}_{\rm m} \cdot \underline{v} \, dx 
	=
	\int_{\Gamma_{\rm N}} \, {\underline{g}}_{\rm m,N}  \cdot \underline{v}   \, dx,
\end{equation}
for all $\underline{v}\in \mathcal{X}^{\rm  u}_{\rm{hf}}$, such that $\underline{v}\cdot \underline{n}|_{\partial \Omega \setminus  \Gamma_{\rm N}}=0$.
\subsection{Finite element formulation}
We resort to an implicit Euler time discretisation scheme, with $J_{\rm{max}}=100$ uniform time steps; the superscript $(\cdot)^{+}$ refers to the new solution (at the current time step $j$, for $j=1,...,J_{\rm{max}}$), while $(\cdot)^{-}$ refers to the solution at the previous time steps:
\begin{equation}
	\left\{
	\begin{array}{l}
		\displaystyle{
			\int_{\Omega} \,  
			2 {\mu} \, \nabla_{\rm s} \, \underline{u}^{+} \, : \,  \nabla_{\rm s} \, \underline{v} \, + \,\left( \lambda \nabla \cdot \underline{u}^{+} \, - \, (2 {\mu} + 3 {\lambda}) \, {\alpha}_{\rm s} T^{+} - b {p}^{+}_{\rm w}  \right) \, \nabla \cdot \underline{v} \, - \, \left({\rho}^0 + m^{+}_{\rm w} \right) \, {\underline{F}}_{\rm m} \cdot \underline{v} \, dx 
		}
		\\[2mm]
		\displaystyle{
			\hspace{3in}
			= \, 
			\int_{\Gamma_{\rm N}} \, \underline{g}^{+}_{\rm m,N}  \cdot \mathbf{v}   \, dx;
		}
		\\[3mm]
		\displaystyle{
			\int_{\Omega} \,
			\frac{1}{\Delta t}(m^{+}_{\rm w} -m^{-}_{\rm{w}})\, \psi \,  + \, 
			{\gamma^{+}} \, (\nabla {p}^{+}_{\rm w}-\rho^{+}_{\rm w}\underline{F}_{\rm m})
			\cdot \nabla \psi \, dx \, = 0;
		}
		\\[3mm]
		\displaystyle{
			\int_{\Omega} \, 
			\Big(
			\left(
			\frac{{h}_{\rm w}}{\Delta t}  (m^{+}_{\rm w}-m^{-}_{\rm{w}}) \, + \, \frac{1}{\Delta t}(\mathcal{Q}^{+}-\mathcal{Q}^{-})
			\, + \,
			\gamma^{+} \,(\nabla {p}^{+}_{\rm w}-\rho^{+}_{\rm w}\underline{F}_{\rm m}) \cdot \underline{F}_{\rm m}
			\right) \, \xi \, - \,
			\left( - \, h^{-}_{\rm w} \, (\nabla p^{+}_{\rm w}-\rho^{+}_{\rm w}{\underline{F}}_{\rm m}) \, + \,{\underline{q}} \right)
			\, \cdot \nabla \xi
		}
		\\[2mm]
		\displaystyle{
			\hspace{3in}
			\,= \, \int_{\Omega} \, \Theta^{+} \, \xi \, dx \,- \, \int_{\partial \Omega} \, \
			g^{+}_{\rm t, N} \, \xi \, dx;
		}
		\\
	\end{array}
	\right.
	\label{eq:FE_THM}
\end{equation}
for all
$\underline{v} \in \mathcal{X}^{\rm  u}_{\rm{hf}}$ such that $\underline{v}\cdot \underline{n}|_{\partial \Omega \setminus \Gamma_N}=0$,
$\psi \in \mathcal{X}^{\rm  p}_{\rm{hf}}$,
$\xi \in \mathcal{X}^{\rm  t}_{\rm{hf}}$,
where 
\begin{equation}
	\label{eq:time_discretization_dependent}
	\left\{
	\begin{array}{l}
		\displaystyle{
			\rho_{\rm w}^+ \, = \, \rho_{\rm w}^- \,
			{\rm exp} \left(\frac{p_{\rm w}^+ - p_{\rm w}^- }{K_{\rm w}} \, - \,
			3 \alpha_{\rm w} (T^+ - T^-)  \right);
		}
		\\[3mm]
		\displaystyle{
			\varphi^+ = b \, - \, (b - \varphi^-) \,
			{\rm exp} \left(
			- (\epsilon_{\rm V}^+ - \epsilon_{\rm V}^-) + 3 \alpha_0 (T^+ - T^-) \,  - \frac{1}{K_s} \, (p_{\rm w}^+ - p_{\rm w}^-) \right);
		}
		\\[3mm]
		\displaystyle{
			h_{\rm w}^+  = 
			h_{\rm w}^- \, + \, 
			C_{\rm w}^{\rm p} \, (T^+  - T^-)  \, + \, 
			\frac{   \beta_{h}^{\rm p} - 3 \alpha_{\rm w} T^+ }{ \rho_{\rm w}^+} 
			\left( p_{\rm w}^+ - p_{\rm w}^- \right);}
		\\[3mm]
		\displaystyle{
			\mathcal{Q}^+ =  \mathcal{Q}^- \,  + \,
			\left( \beta_{\mathcal{Q}}^{\epsilon}  + 3 \alpha_{\rm s} K_0 \, \frac{1}{2}(T^+ + T^-)   \right) \, \left( \epsilon_{\rm V}^+ - \epsilon_{\rm V}^- \right)
			\, -
			\left( \beta_{\mathcal{Q}}^{\rm p}  + 
			3 \alpha_{\rm w,m}^+  \frac{1}{2}(T^+ + T^-)  \right) \, 
			\left( p_{\rm w}^+ \,-\, p_{\rm w}^- \right)
		}
		\\
		\displaystyle{
			\hspace{3in}
			\, + \, C_{\epsilon}^{0,+} \, (T^+ - T^-);}
		\\[3mm]
		\displaystyle{
			m_{\rm w}^+ = \rho_{\rm w}^+ (1 + \epsilon_{\rm V}^+) \, \varphi^+ - \rho_{\rm w}^0 \, \varphi^0.
		}
	\end{array}
	\right.
\end{equation}

\subsection{Choice of the norm}
We equip the FE space $\mathcal{X}_{\rm{hf}}$ with the weighted inner product
\begin{equation}
	(\underline{U}, \underline{U}' ) =\frac{1}{\lambda_{\rm{u}}}\sum_{d=1}^2 (\underline{u}_{d}, \underline{u}'_{d})_{H^{1}(\Omega)}+\frac{1}{\lambda_{\rm{p}}}(p,p')_{H^{1}(\Omega)}+
	\frac{1}{\lambda_{\rm{t}}}(T,T')_{H^{1}(\Omega)},
	\label{eq:norm}
\end{equation}
where the coefficients $\lambda_{\rm{u}}, \lambda_{\rm{p}}, \lambda_{\rm{t}}$ are the largest eigenvalues of the Gramian matrices $\mathbf{C}^{\rm{u}}$, $\mathbf{C}^{\rm{p}}$, $\mathbf{C}^{\rm{t}}$ associated to displacement, pressure and temperature, respectively. Similarly to \cite{taddei2019offline}, the inner product \eqref{eq:norm} is  motivated by the need for properly taking into account the contributions of displacement, pressure and temperature, which are characterised by different magnitudes and different units.

\subsection{Parametrization}
\label{sec:THM_parametrization}
We consider a vector of four parameters: the Young's modulus $E$ and the Poisson's ratio $\nu$ in the region UA, the thermic factor $\tau$ and the constant $C_{\rm{al}}$ in \eqref{eq:gtN}. For all parameters, we define the parameter domain $\mathcal{P}$ by considering variations of $\pm 15 \%$ with respect to the nominal value reported in Table \ref{table:further}. % and \eqref{eq:gtN}.

\section{Numerical results}
\label{sec:numerics}

We measure performance through the discrete $L^2(0,T_{\rm{f}};\mathcal{X}_{\rm{hf}})$ relative error
\begin{equation}
	E_{\mu}:=\frac{\sqrt{\displaystyle{\sum_{j=1}^{J_{\rm{max}}}}{
				\left(t^{(j)}-t^{(j-1)}\right) \big\| \underline{U}^{(j)}_{\rm{hf},\mu}-\widehat{\underline{U}}^{(j)}_{\mu} \big\|^2}}}{\sqrt{\displaystyle{\sum_{j=1}^{J_{\rm{max}}}}{
				\left(t^{(j)}-t^{(j-1)}\right) \| \underline{U}^{(j)}_{\rm{hf},\mu}\|^2}}}
			\label{eq:relative_L2error}
\end{equation}
for any $\mu \in \mathcal{P}$.
Similarly, we denote by $E^{\rm{u}}_{\mu}$, $E^{\rm{p}}_{\mu}$ and $E^{\rm{t}}_{\mu}$ the discrete relative $L^2(0,T_{\rm{f}};\mathcal{X}_{\rm{hf}})$ errors associated with the estimate of displacement, pressure and temperature, respectively.

\subsection{Solution reproduction problem}
\label{sec:solution_rep_pb}
We first present numerical results for a fixed configuration of parameters $\bar{\mu} \in \mathcal{P}$ to validate the ROM described in section \ref{sec:method}.  We consider $\bar{\mu}$ equal to the centroid of $\mathcal{P}$. We perform data compression based on the whole set of snapshots, i.e. $|I_{\rm{s}}|=J_{\rm{max}}=100$.

\subsubsection{Data compression: POD}
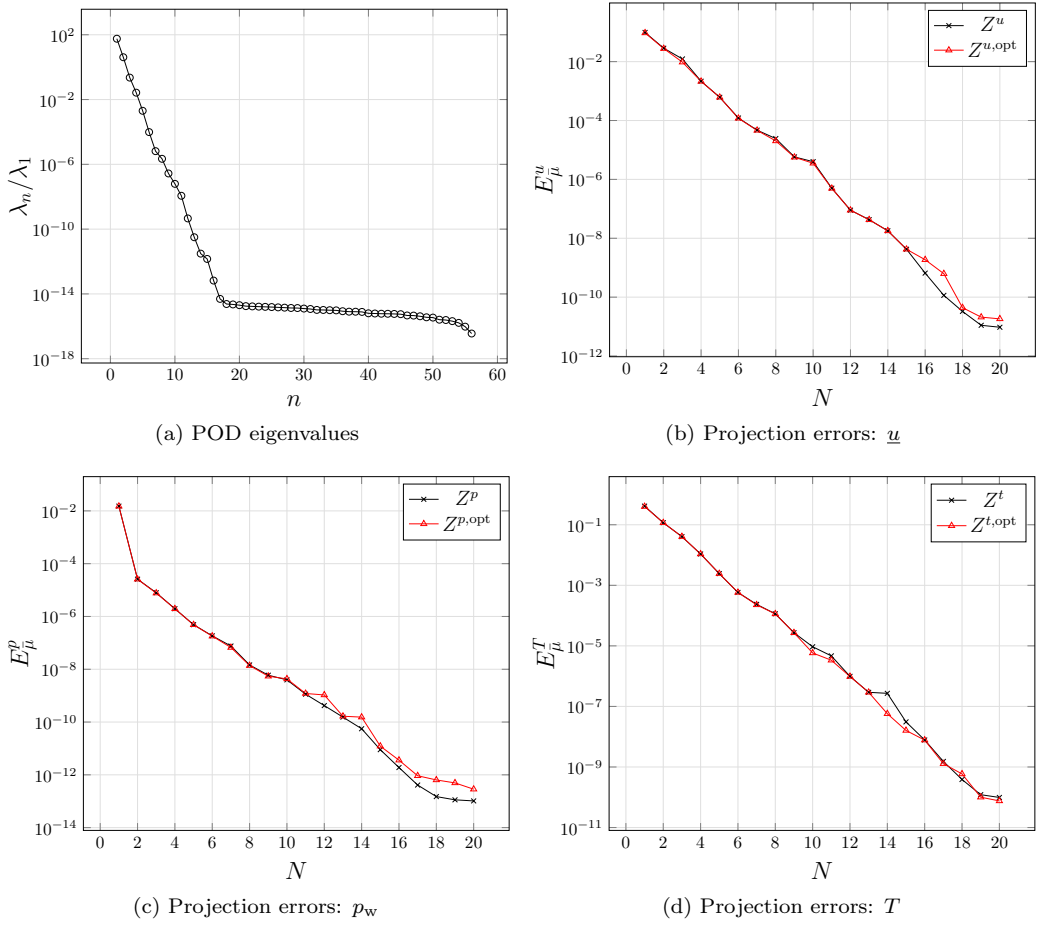
\begin{figure} %[H]
	\centering
	\subfloat[POD eigenvalues]{
		\label{fig:eigenvalues}
	\begin{tikzpicture}[scale=.65]
		\begin{loglogaxis}[
			xmode=linear,
			ymode=log,
			xlabel={\Large{$n$}},
			ylabel={\Large{$\lambda_n/\lambda_1$}},
			width=0.6\textwidth,
			grid=both,
			minor grid style={gray!25},
			major grid style={gray!25},
			]
			]
			\addplot[ % <-- plot options
			black,
			mark=o,
			mark options={ color = black},
			]  %
			table{dat/srp/eigenvalues.dat};
		\end{loglogaxis}
	\end{tikzpicture}
}
\subfloat[Projection errors: $\underline{u}$]{
	\label{fig:proj-u}
		\begin{tikzpicture}[scale=.65]
			\begin{loglogaxis}[
				xmode=linear,
				ymode=log,
				xlabel={\Large{$N$}},
				ylabel={\Large{$E^{u}_{\bar{\mu}}$}},
				width=0.6\textwidth,
				grid=both,
				minor grid style={gray!25},
				major grid style={gray!25},
				]
				]
				\addplot[ % <-- plot options
				black,
				mark=x,
				mark options={ color = black},
				]  %
				table{dat/srp/POD/e_u.dat};
				\addlegendentry{\large{$Z^{u}$}};
				\addplot[ % <-- plot options
				red,
				mark=triangle,
				mark options={ color = red},
				]  %
				table{dat/srp/POD/e_u_opt.dat};
				\addlegendentry{\large{$Z^{u, \rm{opt}}$}};
			\end{loglogaxis}
		\end{tikzpicture}
	}
\quad 
\subfloat[Projection errors: $p_{\rm{w}}$]{
	\label{fig:proj-p}
		\begin{tikzpicture}[scale=.65]
			\begin{loglogaxis}[
				xmode=linear,
				ymode=log,
				xlabel={\Large{$N$}},
				ylabel={\Large{$E^{p}_{\bar{\mu}}$}},
				width=0.6\textwidth,
				grid=both,
				minor grid style={gray!25},
				major grid style={gray!25},
				]
				]
				\addplot[ % <-- plot options
				black,
				mark=x,
				mark options={ color = black},
				]  %
				table{dat/srp/POD/e_p.dat};
				\addlegendentry{\large{$Z^{p}$}};
				\addplot[ % <-- plot options
				red,
				mark=triangle,
				mark options={ color = red},
				]  %
				table{dat/srp/POD/e_p_opt.dat};
				\addlegendentry{\large{$Z^{p,\rm{opt}}$}};
			\end{loglogaxis}
		\end{tikzpicture}
	}
\subfloat[Projection errors: $T$]{
	\label{fig:proj-T}
		\begin{tikzpicture}[scale=.65]
			\begin{loglogaxis}[
				xmode=linear,
				ymode=log,
				xlabel={\Large{$N$}},
				ylabel={\Large{$E^{T}_{\bar{\mu}}$}},
				width=0.6\textwidth,
				grid=both,
				minor grid style={gray!25},
				major grid style={gray!25},
				]
				]
				\addplot[ % <-- plot options
				black,
				mark=x,
				mark options={ color = black},
				]  %
				table{dat/srp/POD/e_T.dat};
				\addlegendentry{\large{$Z^{t}$}};
				\addplot[ % <-- plot options
				red,
				mark=triangle,
				mark options={ color = red},
				]  %
				table{dat/srp/POD/e_T_opt.dat};
				\addlegendentry{\large{$Z^{t,\rm{opt}}$}};
			\end{loglogaxis}
		\end{tikzpicture}
	}
	\caption{\protect\subref{fig:eigenvalues}: exponential decay of POD eigenvalues. \protect \subref{fig:proj-u}, \protect\subref{fig:proj-p}, \protect\subref{fig:proj-T}: projection errors computed through \eqref{eq:global_POD} (in black) and \eqref{eq:optimal_POD_u}-\eqref{eq:optimal_POD_T}(in red) for  increasing numbers of POD modes.}
	\label{fig:projections}
	\end{figure}

In Figure \ref{fig:projections} we compare performance of the global POD based on the weighted inner product $(\cdot, \cdot)$ with the performance of the component-wise POD. More precisely, we define $\underline{Z}$ such that 
\begin{equation}
	[\underline{Z},\boldsymbol{\lambda}]=\texttt{POD} \left(  
	\{    \underline{U}_{\rm hf, \bar{\mu}}^{(j)}  \}_{j\in I_{\rm{s}}},
	(\cdot,\cdot), tol_{\rm pod}  \right),
	\label{eq:global_POD}
\end{equation}
 and we  extract the displacement, pressure and temperature components $\underline{Z}^{\rm{u}}$, $\underline{Z}^{\rm{p}}$, $\underline{Z}^{\rm{t}}$. Then , we denote the "optimal" (in a discrete $L^2$ sense) spaces 
 \begin{align}
 	&[\underline{Z}^{u,\rm{opt}},\boldsymbol{\lambda}^{u,\rm{opt}}]=\texttt{POD}
 	\left( \{    \underline{u}_{\rm hf, \bar{\mu}}^{(j)}  \}_{j \in \texttt{I}_{\rm{s}}}, 
 	(\cdot,\cdot)_{H^1}, tol_{\rm pod}  \right);
 	\label{eq:optimal_POD_u}\\
 	&[\underline{Z}^{p,\rm{opt}},\boldsymbol{\lambda}^{p,\rm{opt}}]=\texttt{POD}
 	\left( \{   p_{\rm hf, \bar{\mu}}^{(j)}  \}_{j \in \texttt{I}_{\rm{s}}}, 
 	(\cdot,\cdot)_{H^1}, tol_{\rm pod}  \right);\\
 	&[\underline{Z}^{T,\rm{opt}},\boldsymbol{\lambda}^{T,\rm{opt}}]=\texttt{POD}
 	\left( \{   T_{\rm hf, \bar{\mu}}^{(j)}  \}_{j \in \texttt{I}_{\rm{s}}}, 
 	(\cdot,\cdot)_{H^1}, tol_{\rm pod}  \right).
 	\label{eq:optimal_POD_T}
 \end{align}
 In Figure \ref{fig:projections} \subref{fig:eigenvalues} we show the behaviour of the POD eigenvalues in \eqref{eq:global_POD}; in Figure \ref{fig:projections}\subref{fig:proj-u}, \subref{fig:proj-p}, \subref{fig:proj-T} we compare the relative projection errors associated with $\underline{Z}^{\rm{u}}$ and $\underline{Z}^{\rm{u, opt}}$, $\underline{Z}^{\rm{p}}$, $\underline{Z}^{\rm{p, opt}}$ and $\underline{Z}^{\rm{t}}$ and $\underline{Z}^{\rm{t, opt}}$.
We observe that the projection errors are nearly the same for all the three state variables: this obervation suggests to consider a single reduced space to approximate the solution field.

\subsubsection{Hyper-reduction}

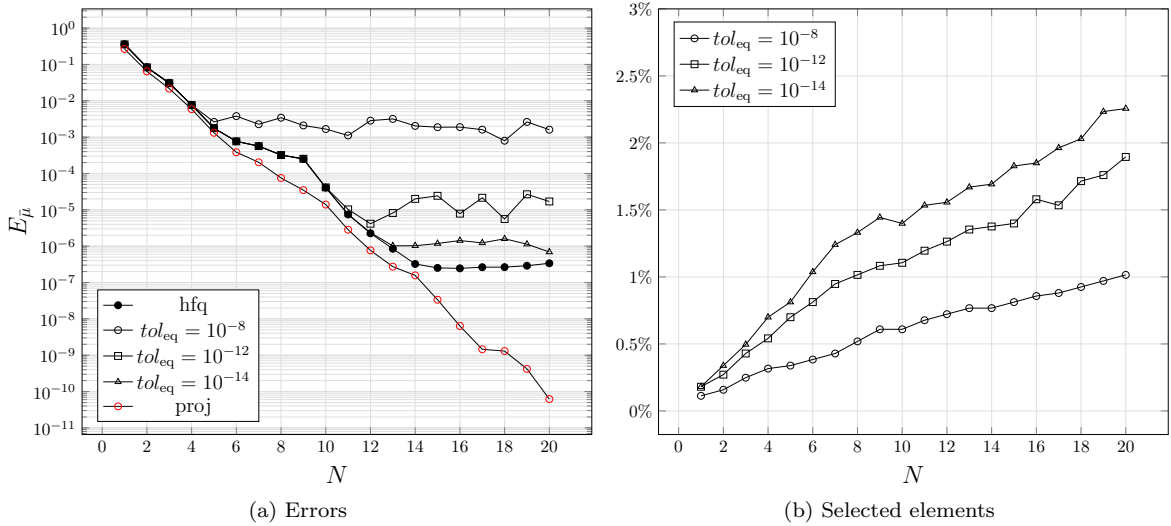
\begin{figure} %[H]
	\centering
	\subfloat[Errors]{	
		\label{fig:hyper_errors}
		\begin{tikzpicture}[scale=.65]
			\begin{loglogaxis}[
				xmode=linear,
				ymode=log,
				xlabel={\Large{$N$}},
				ylabel={\Large{$E_{\bar{\mu}}$}},
				title={},
				grid=both,
				minor grid style={gray!25},
				major grid style={gray!25},
				width=0.7\linewidth,
				legend pos=south west,
				]
				]
				\addplot[ % <-- plot options
				mark=*,
				mark options={color = black}
				]  %
				table{dat/srp/error_hfq.dat};
				\addlegendentry{\large{hfq}};
				\addplot[ % <-- plot options
				mark=o,
				mark options={color =black}
				]  %
				table{dat/srp/error_eq_tol1e-8.dat};
				\addlegendentry{\large{$tol_{\rm{eq}}=10^{-8}$}};
				\addplot[ % <-- plot options
				mark=square,
				mark options={color = black}
				]  %
				table{dat/srp/error_eq_tol1e-12.dat};
				\addlegendentry{\large{$tol_{\rm{eq}}=10^{-12}$}};
				\addplot[ % <-- plot options
				mark=triangle,
				mark options={color = black}
				]  %
				table{dat/srp/error_eq_tol1e-14.dat};
				\addlegendentry{\large{$tol_{\rm{eq}}=10^{-14}$}};
					\addplot[ % <-- plot options
				mark=o,
				mark options={color = red}
				]  %
				table{dat/srp/eproj.dat};
				\addlegendentry{\large{proj}};
			\end{loglogaxis}
		\end{tikzpicture}
	}
	\subfloat[Selected elements]{
		\label{fig:elements}
	\begin{tikzpicture}[scale=.65]
		\begin{axis}[
			xlabel={\Large{$N$}},
			ylabel={},
			title={},
			grid=both,
			minor grid style={gray!25},
			major grid style={gray!25},
			width=0.7\linewidth,
			yticklabel={\pgfmathparse{\tick}\pgfmathprintnumber{\pgfmathresult}\%},
			ymax=3.0,
			legend cell align={left},
			legend pos=north west,
			]
			]
			\addplot[ % <-- plot options
			mark=o,
			mark options={color = black}
			]  %
			table{dat/srp/elements_tol1e-8.dat};
				\addlegendentry{\large{$tol_{\rm{eq}}=10^{-8}$}};
				\addplot[ % <-- plot options
			mark=square,
			mark options={color =black}
			]  %
			table{dat/srp/elements_tol1e-12.dat};
				\addlegendentry{\large{$tol_{\rm{eq}}=10^{-12}$}};
				\addplot[ % <-- plot options
			mark=triangle,
			mark options={color = black}
			]  %
			table{dat/srp/elements_tol1e-14.dat};
			\addlegendentry{\large{$tol_{\rm{eq}}=10^{-14}$}};
		%	\legend{$tol_{EQ}=10^{-8}$,$tol_{EQ}=10^{-12}$,$tol_{EQ}=10^{-14}$.}
		\end{axis}
	\end{tikzpicture}
}
	\caption{solution reproduction problem. \protect\subref{fig:hyper_errors}: errors associated to projection error (proj), Galerkin with high-fidelity quadrature (hfq) and Galerkin with empirical quadrature for several choices of $tol_{\rm{eq}}$ with respect to the ROM dimension $N$. \protect\subref{fig:elements}: percentage of selected elements for several $tol_{\rm{eq}}$.}
	\label{fig:hyper-reduction}
\end{figure}

In Figure \ref{fig:hyper-reduction}\subref{fig:hyper_errors} we show the performance of the Galerkin ROM with and without hyper-reduction. We distinguish between the high-fidelity quadrature rule, abbreviated as $\rm{hfq}$, and the empirical quadrature rule for several tolerances $tol_{\rm{eq}}$. 
We also add as a reference, the relative projection error.
Figure \ref{fig:hyper-reduction}\subref{fig:elements} shows the percentage of selected elements $\frac{\rm{Q}}{N_{\rm{e}}}\times 100 \%$  for the same choices of the tolerance $tol_{\rm{eq}}$.
We observe that the empirical quadrature procedure is able to significantly reduce the size  of the mesh used for online calculations without compromising accuracy. The plateau for $N \gtrsim 14$ is due to the tolerance of the Newton iterative solver.

\begin{figure}%[H]
	\centering
	\subfloat[$tol_{\rm{eq}}=10^{-14}$]{
		\label{fig:tol1e-14}
		\includegraphics[width=.5\textwidth]{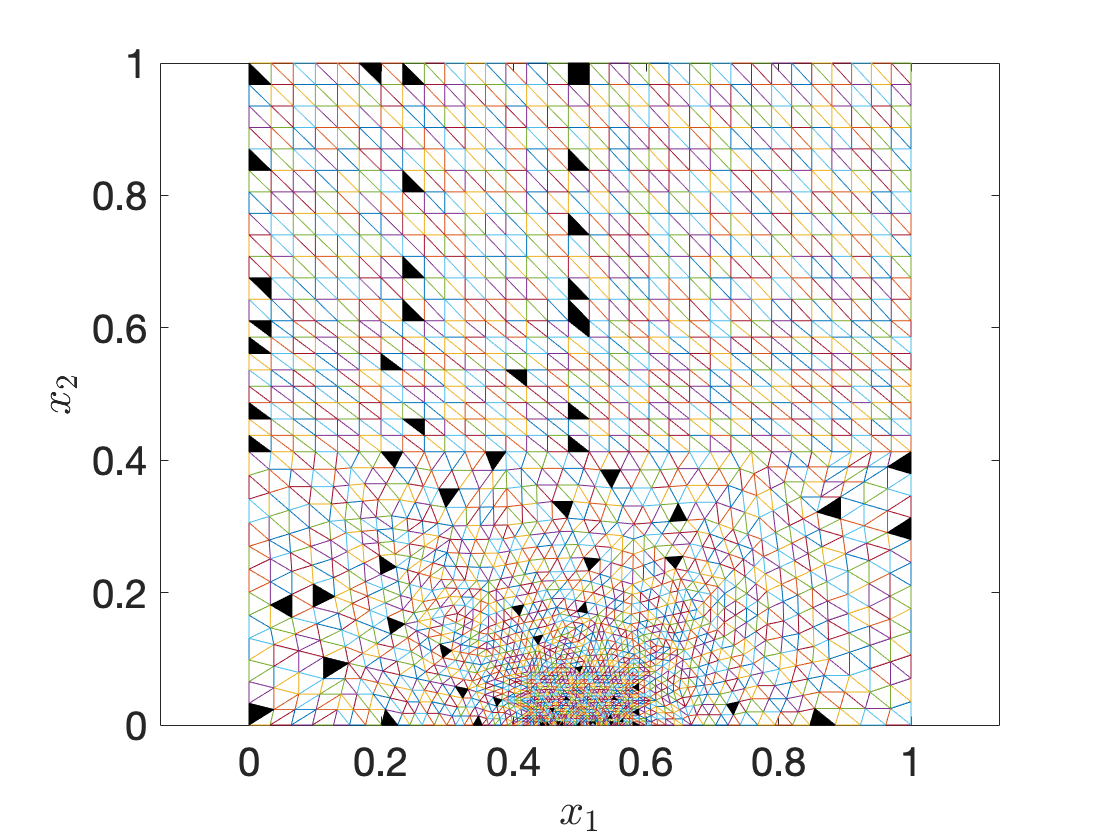}
	}
	\subfloat[$tol_{\rm{eq}}=10^{-10}$]{
		\label{fig:tol1e-10}
		\includegraphics[width=.5\textwidth]{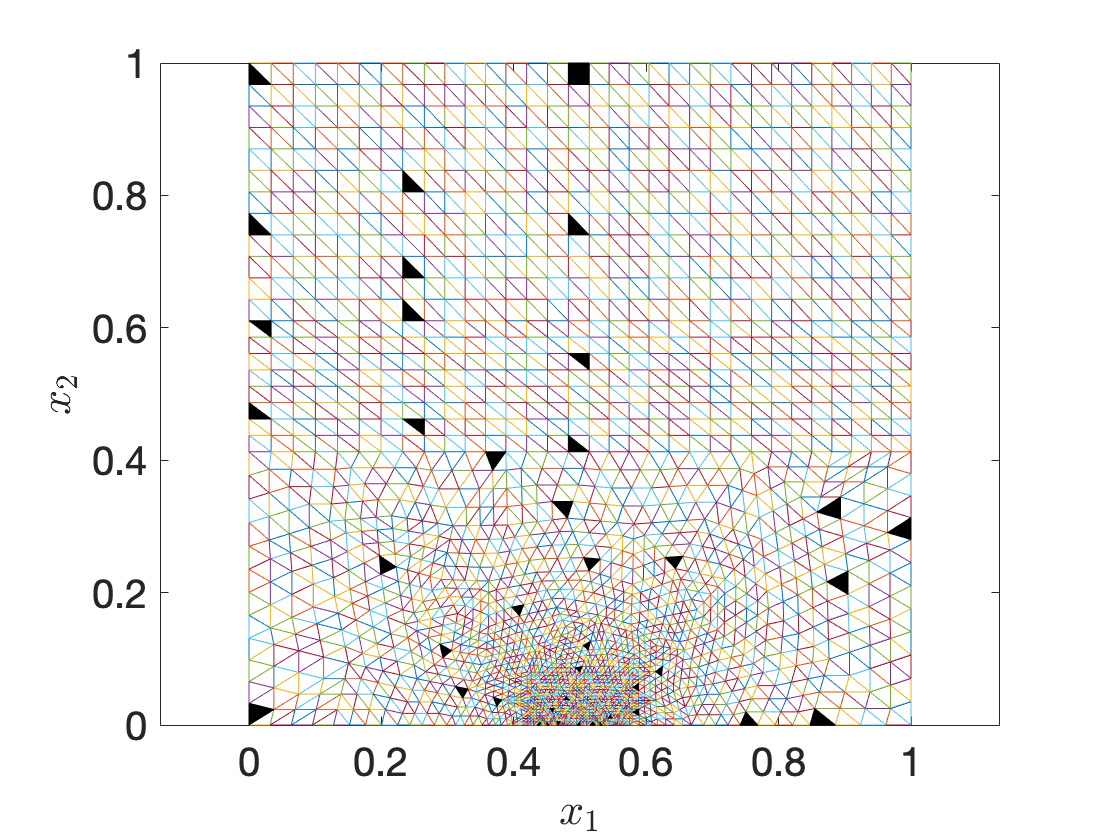}
	}
	\caption{solution reproduction problem. Reduced mesh for two choices of the empirical quadrature tolerance.}
	\label{fig:selected_elem}
\end{figure}

In Figure \ref{fig:selected_elem}, we show the selected grid elements for two choices of the EQ tolerance  value $tol_{\rm{eq}}$ and for $N=12$. We observe that the sampled elements are distributed over the whole domain with a slight prevalence of elements in the proximity of the alveoli.

\subsection{Parametric problem}
\label{sec:parametric pb}
We present results for the parametric case. 
We denote by $\Xi_{\rm{train}} \subset \mathcal{P}$ the training set used to build the ROM and by $\Xi_{\rm test} \subset \mathcal{P}$  the test set used to assess performance.
Both sets consist of independent identically distributed samples of a uniform distribution in $\mathcal{P}$, with $|\Xi_{\rm{train}}|=n_{\rm{train}}=50$ and $|\Xi_{\rm test}|=n_{\rm{test}}=10$. We also set $tol_{\rm{POD}}=10^{-7}$ in \eqref{eq:POD_maxN} and in \eqref{eq:HAPOD_N} for data compression, and we set $tol_{\rm POD,res}=10^{-5}$ in \eqref{eq:POD_maxN} for the construction of the empirical test space. We set $I_{\rm s} \subset \{1,...,J_{\rm{max}}\}$ with $|I_{\rm s}|=20$. EQ rules are depicted usign the tolerance $tol_{\rm{eq}}=10^{-12}$ (cf . Algorithm \ref{alg:ROM_construction}).

%The reduced bases are found through a POD with a fixed tolerance $tol_{\rm{POD}}=10^{-7}$. The POD tolerance used to define the reduced test space is $10^{-5}$. \\
%The empirical quadrature tolerance is $tol_{\rm{eq}}=10^{-12}$ for the reduced order model used both to compute the reduced solutions and the error indicator. 

\subsubsection{Error estimation}

\begin{figure}[h!]
	\centering
	\subfloat[H-POD]{
	\label{fig:HPOD_correlation}
	\begin{tikzpicture}[scale=.6]
		\begin{loglogaxis}[
				xmode=log,
				ymode=log,
			xlabel={\LARGE{$\Delta_{\mu}$}},
			ylabel={\LARGE{$E_{\mu}$}},
			title={},
			grid=both,
			minor grid style={gray!25},
			major grid style={gray!25},
			width=0.7\linewidth,
			ymin=1e-6,ymax=1,
			xmin=1e-6,xmax=1e-2,
			]
			]
			\addplot[ % <-- plot options
			only marks,
			mark=o,
			mark options={color = black}
			]  %
			table{dat/par_pb/DeltaErr.dat};
			%\addlegendentry{Dual residual,EQ error};
		\end{loglogaxis}
	\end{tikzpicture}
}
	\subfloat[HA-POD]{
	\label{fig:HAPOD_correlation}
	\begin{tikzpicture}[scale=.6]
		\begin{loglogaxis}[
			xmode=log,
			ymode=log,
			xlabel={\LARGE{$\Delta_{\mu}$}},
			ylabel={\LARGE{$E_{\mu}$}},
			title={},
			grid=both,
			minor grid style={gray!25},
			major grid style={gray!25},
			width=0.7\linewidth,
			ymin=1e-6,ymax=1,
			xmin=1e-6,xmax=1e-2,
			]
			]
			\addplot[ % <-- plot options
			only marks,
			mark=o,
			mark options={color = black}
			]  %
			table{dat/par_pb/hapod/DeltaErr.dat};
			%\addlegendentry{Dual residual,EQ error};
		\end{loglogaxis}
	\end{tikzpicture}
}
\caption{parametric problem: correlation between the time-average residual indicator \eqref{eq:error_indicator}  and true relative  errors \eqref{eq:relative_L2error}.}
\label{fig:correlation}
\end{figure}
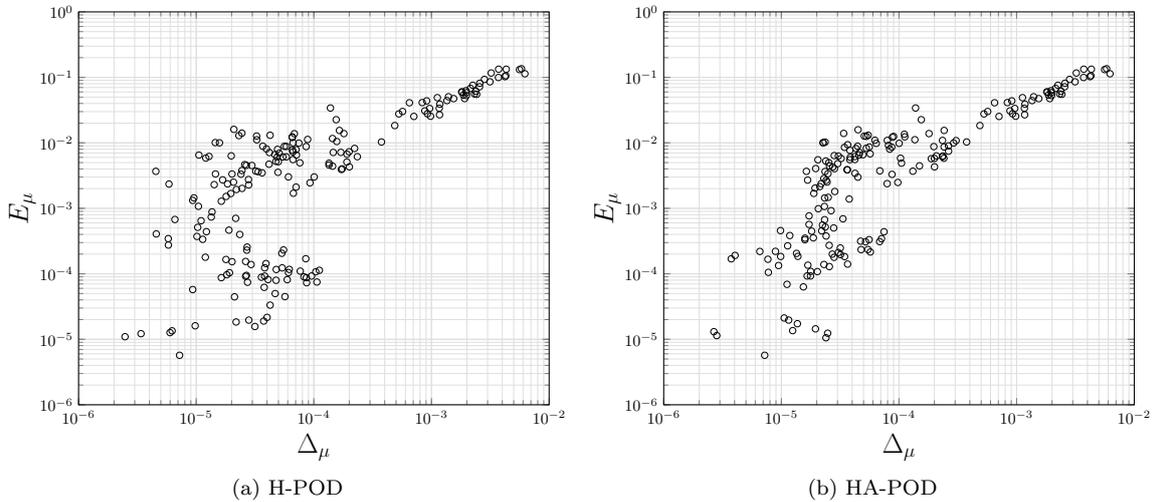

In Figure \ref{fig:correlation} 
we compare the dual residual and several EQ errors for each parameter $\mu$ in the training set $\Xi_{\rm train}$  and for different dimensions of the reduced space that is progressively updated during the execution of the POD-Greedy algorithm. In particular, we show results in two cases: the hierarchical POD-Greedy (H-POD) and the hierarchical approximate POD-Greedy (denoted as HA-POD). 
Figures \ref{fig:correlation}\subref{fig:HPOD_correlation} and \ref{fig:correlation}\subref{fig:HAPOD_correlation}  show for both H-POD and HA-POD  to what extent  the residual-based error indicator defined in \eqref{eq:error_indicator} is correlated with  the relative error \eqref{eq:relative_L2error}.
We observe that   for values of the indicators that are larger than $10^{-3}$,   correlation is very high, while for smaller values   correlation is much weaker. 

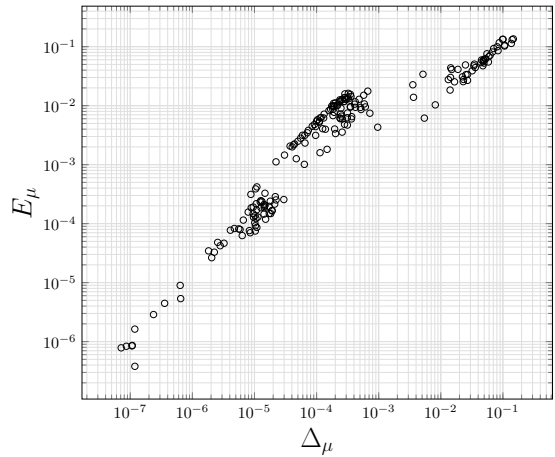
\begin{figure}[h!]
	\centering
\begin{tikzpicture}[scale=.6]
			\begin{loglogaxis}[
				xmode=log,
				ymode=log,
				xlabel={\LARGE{$\Delta_{\mu}$}},
				ylabel={\LARGE{$E_{\mu}$}},
				title={},
				grid=both,
				minor grid style={gray!25},
				major grid style={gray!25},
				width=0.7\linewidth,
				%ymin=1e-6,ymax=1,
				%xmin=1e-6,xmax=1e-2,
				]
				]
				\addplot[ % <-- plot options
				only marks,
				mark=o,
				mark options={color = black}
				]  %
				table{dat/par_pb/strong_PODGreedy/DeltaErr.dat};
				%\addlegendentry{Dual residual,EQ error};
			\end{loglogaxis}
		\end{tikzpicture}
		\caption{parametric problem: correlation between residual indicator \eqref{eq:L2_residual} and true relative errors \ref{eq:relative_L2error}.}
	\label{fig:correlation_strongGreedy}
\end{figure}

To provide a a concrete reference, 
in Figure  \ref{fig:correlation_strongGreedy}
we   investigate the correlation between the relative error \eqref{eq:relative_L2error} and the time-discrete $L^2(0, T_{\rm{f}}; \mathcal{X}'_{\rm{hf},0})$  residual indicator defined in   \eqref{eq:L2_residual}: we observe that the indicator in 
 \eqref{eq:L2_residual} is significantly more accurate, particularly for small values of the error. As stated in section \ref{sec:method}, the residual indicator  \eqref{eq:L2_residual} is considerably more expensive  in terms of both memory and computational costs.

 %A more precise but expensive error indicator is taken into account to complete the numerical results: in Figure\ref{fig:residual_indicator} we compute 

%Comparing the two graphs in \ref{fig:correlation}\subref{fig:HPOD_correlation} and \ref{fig:correlation}\subref{fig:HAPOD_correlation}, we can observe that

%We can also say that for small values of the error indicator (below the order of $10^{-4}$), in the HA-POD case the correlation seems to be slightly higher than in the H-POD case.

\subsubsection{POD-Greedy sampling}

\begin{figure}[h]
		\centering
		\subfloat[Iteration $\rm{it}=1$; $N=15$, $Q=74$, $Q_{\rm r}=16$]{
		\includegraphics[width=0.4\textwidth]{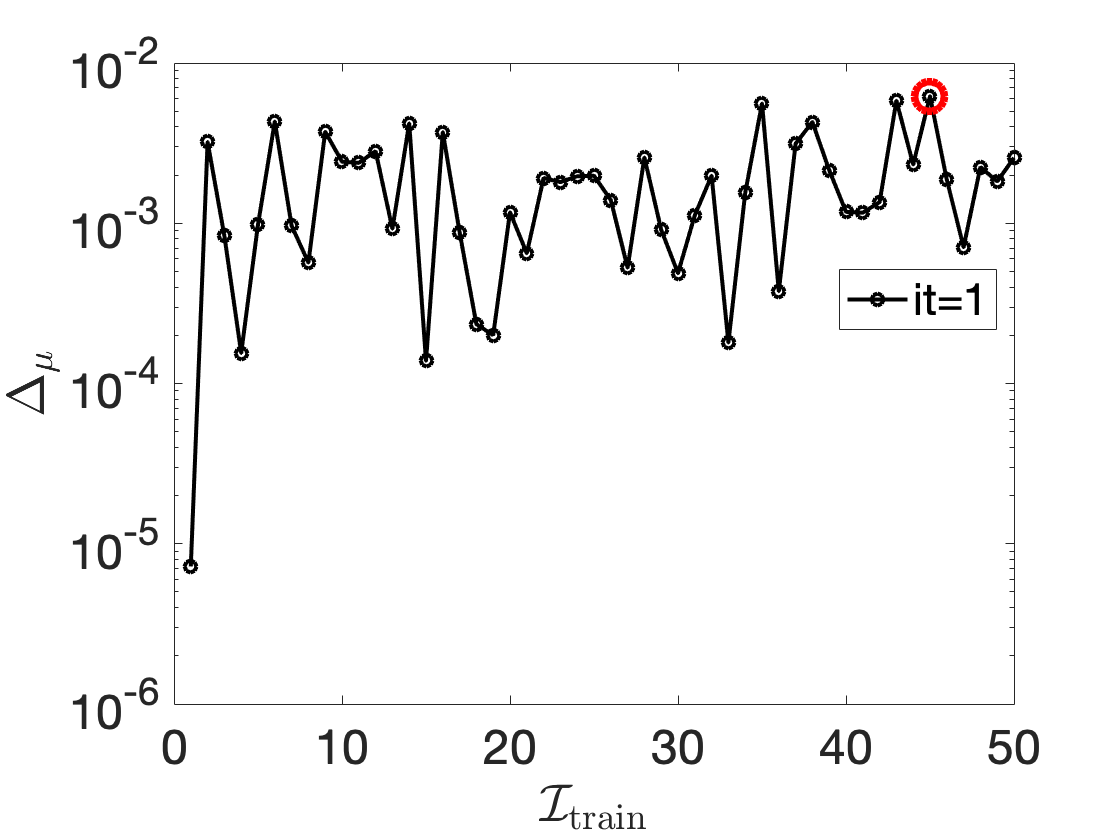}
	}
\subfloat[Iteration $\rm{it}=2$; $N=26$, $Q=123$, $Q_{\rm r}=18$]{
		\includegraphics[width=0.4\textwidth]{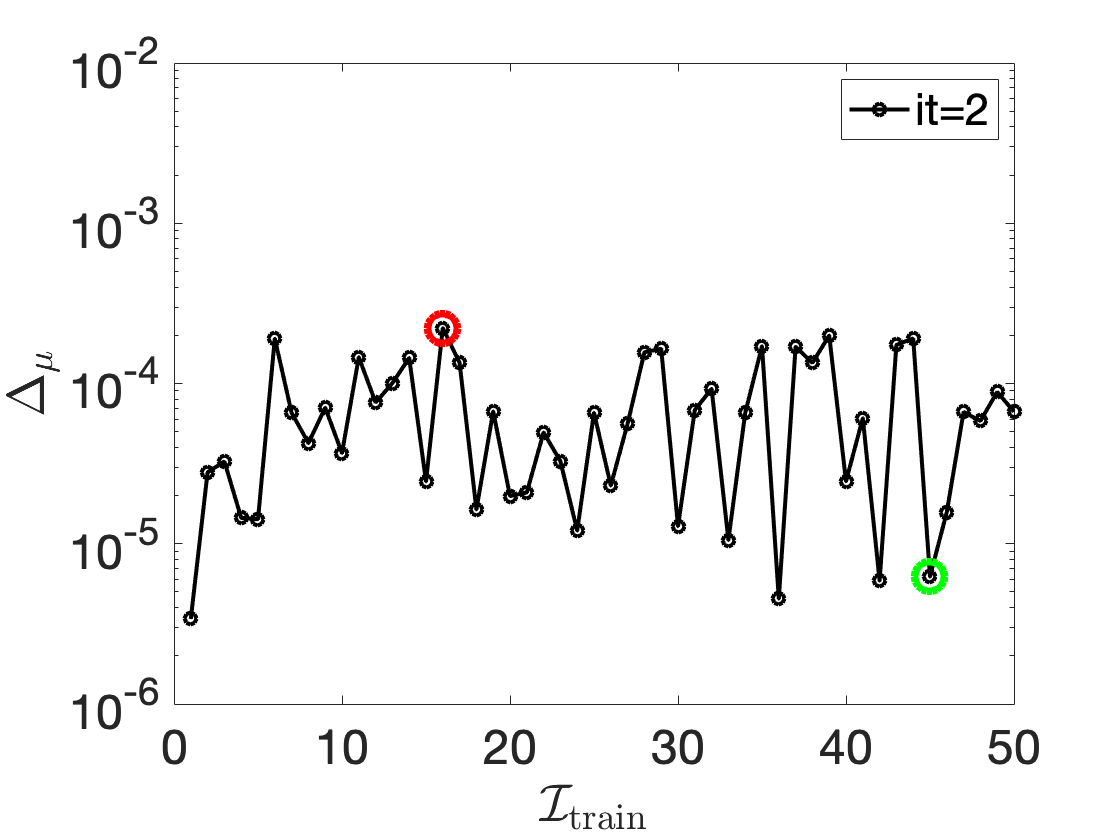}
	}
\quad 
\subfloat[Iteration $\rm{it}=3$; $N=35$, $Q=155$, $Q_{\rm r}=22$]{
		\includegraphics[width=0.4\textwidth]{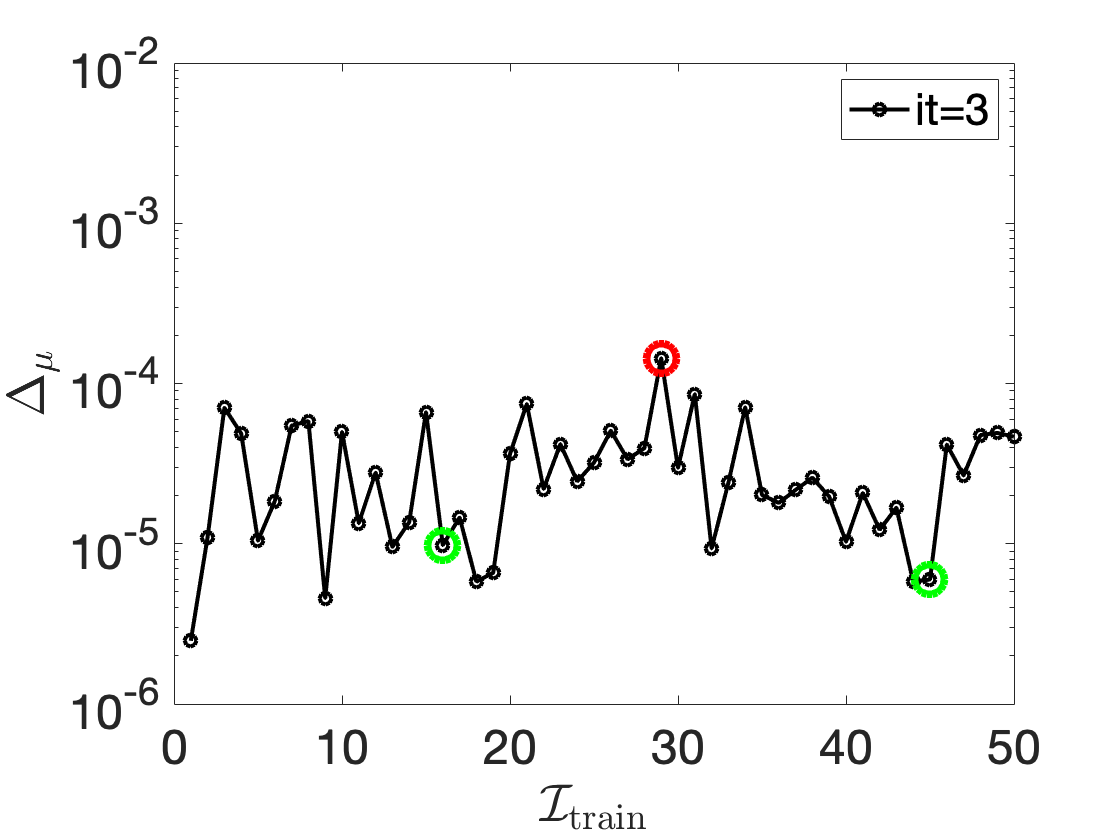}
	}
\subfloat[Iteration $\rm{it}=4$; $N=43$, $Q=169$, $Q_{\rm r}=18$]{
		\includegraphics[width=0.4\textwidth]{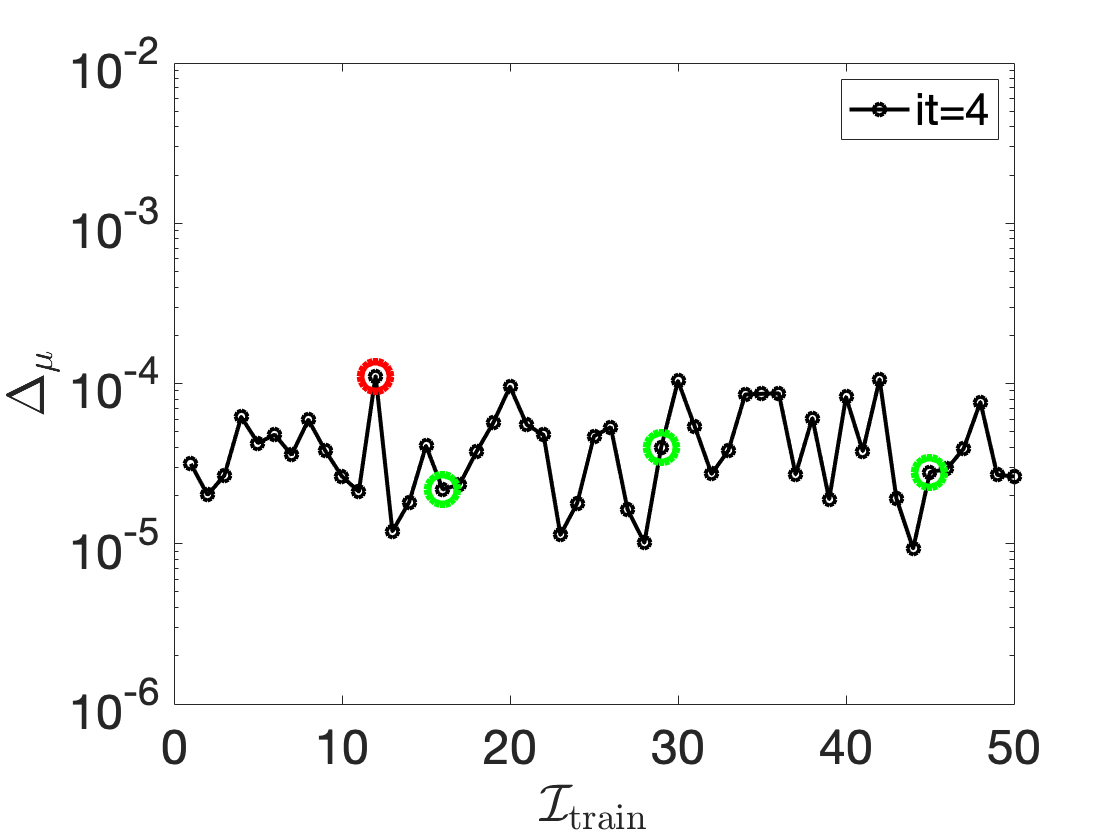}
	}
		\caption{parametric problem: POD-Greedy algorithm convergence history in the H-POD case.}
		\label{fig:H-PODgreedy}
	\end{figure}

\begin{figure}[h]
	\centering
	\subfloat[$\rm{it}=1$; $N=15$,$Q=74$,$Q_{\rm r}=16$]{
		\includegraphics[width=0.4\textwidth]{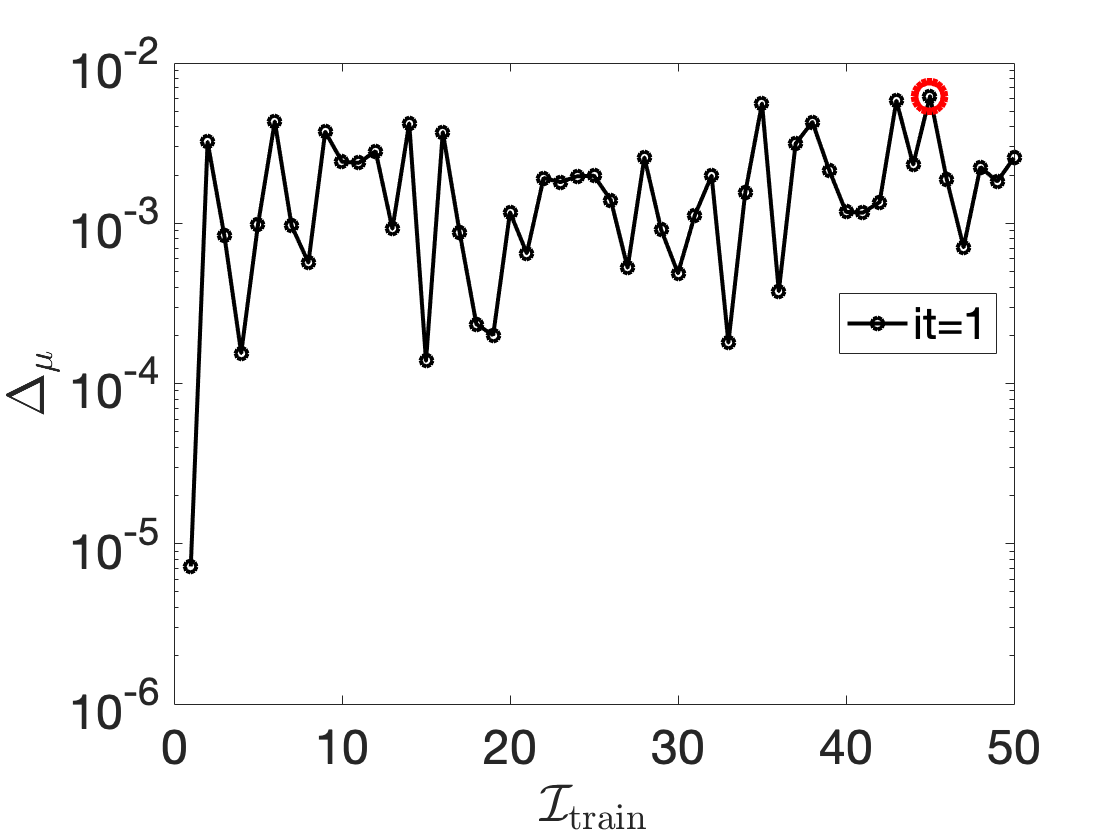}
	}
	\subfloat[$\rm{it}=2$; $N=25$,$Q=120$, $Q_{\rm r}=19$]{
		\includegraphics[width=0.4\textwidth]{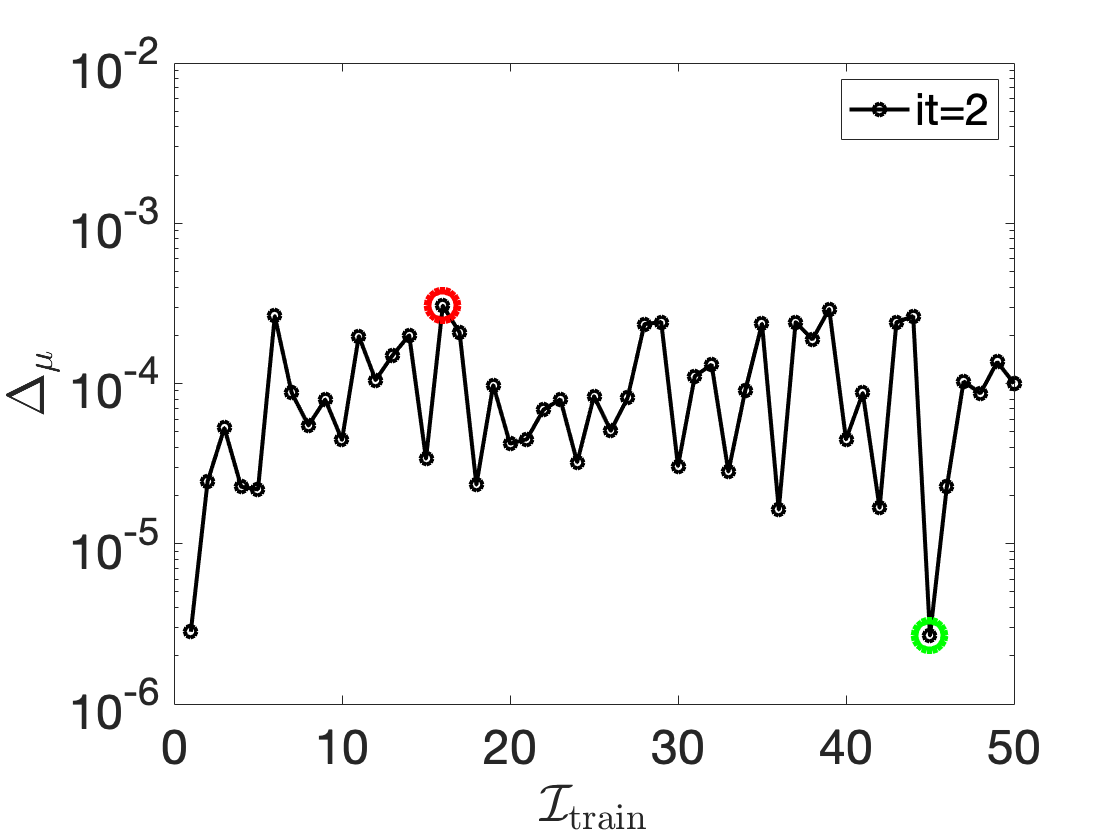}
	}
	\quad 
	\subfloat[$\rm{it}=3$; $N=31$,$Q=135$, $Q_{\rm r}=21$]{
		\includegraphics[width=0.4\textwidth]{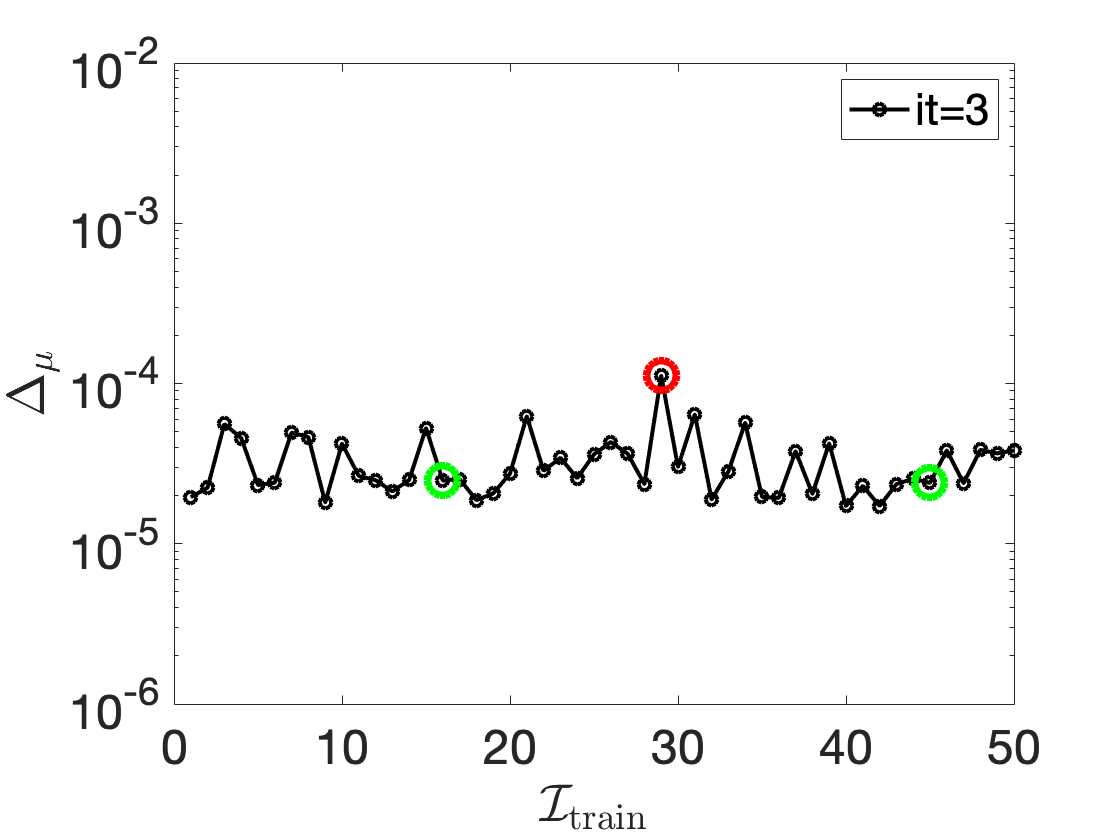}
	}
	\subfloat[Iteration $\rm{it}=4$; $N=37$,$Q=156$, $Q_{\rm r}=19$]{
		\includegraphics[width=0.4\textwidth]{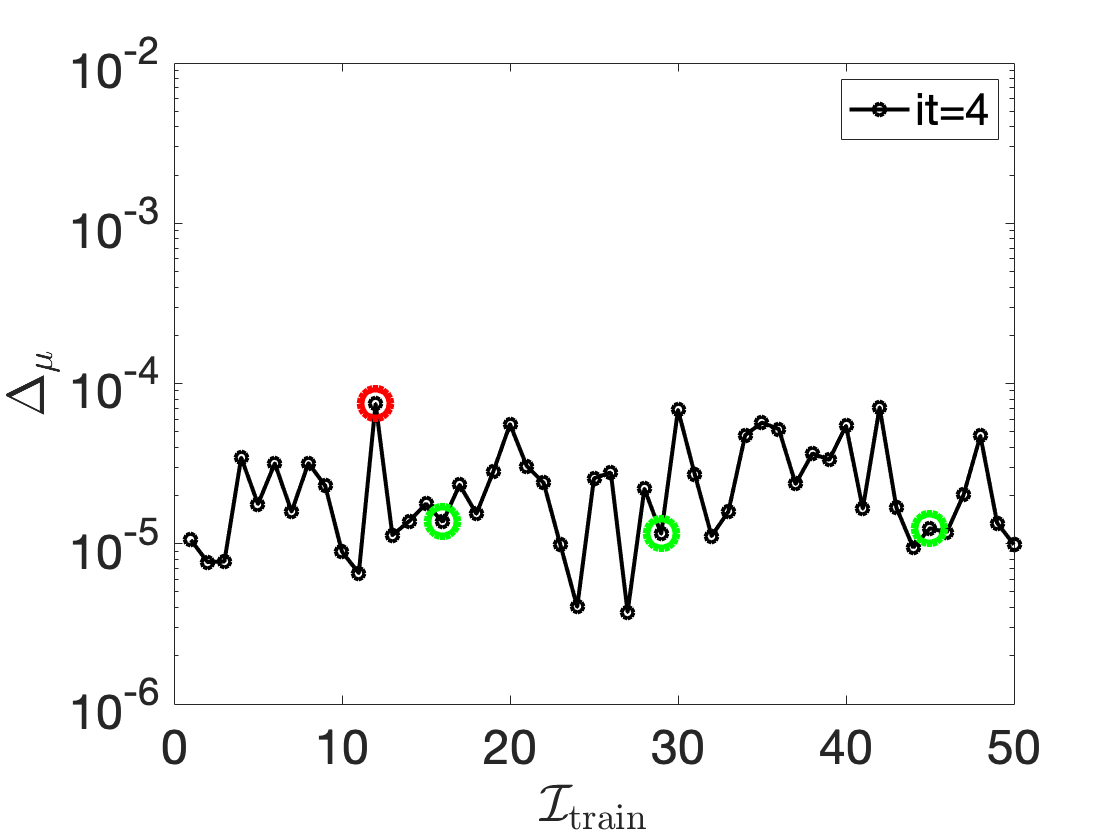}
	}
	\caption{parametric problem: POD-Greedy algorithm convergence history in the HA-POD case. }
	\label{fig:HA-PODgreedy}
\end{figure}

In Figures \ref{fig:H-PODgreedy} and \ref{fig:HA-PODgreedy} we show the POD-Greedy algorithm convergence history, for both the hierarchical and approximate hierarchical PODs. At each iteration of the  algorithm, until convergence, the error indicator $\Delta_{\mu}$ is illustrated with respect to training parameter indices $\mathcal{I}_{\rm{train}}=\{1, ..., |\Xi_{\rm{train}}|\}$ .
 At each iteration the selected parameter $\mu^{\star}$ is marked in red, while the previously selected parameters are marked in green. We also report the dimension of the updated reduced space and the number of sampled elements.

\subsubsection{Predictive tests}

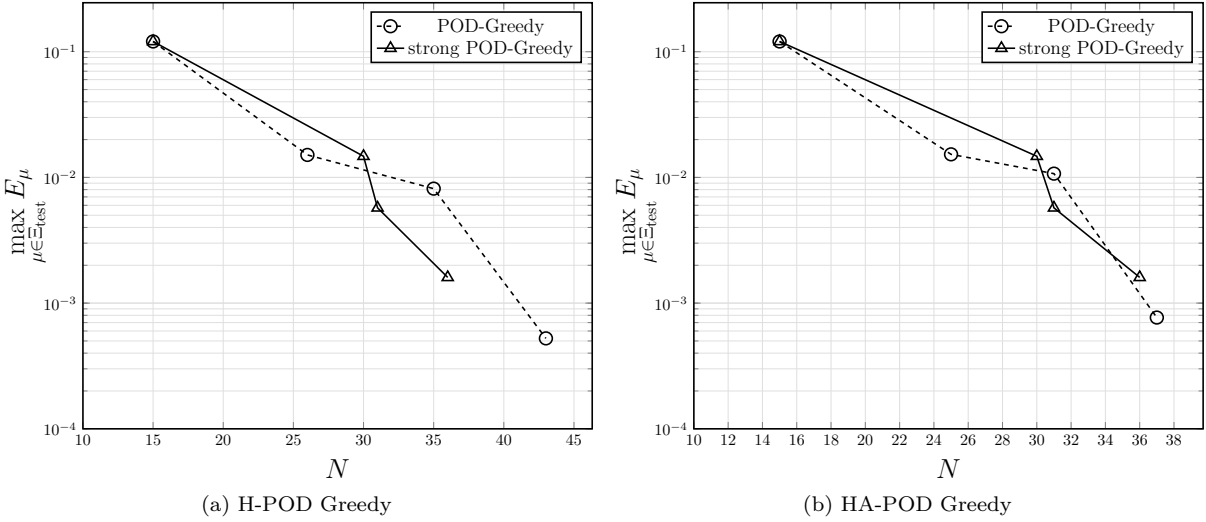
\begin{figure}[h!]
\subfloat[H-POD Greedy]{	
	\label{fig:HPOD_predictions}
	\begin{tikzpicture}[scale=.6]
		\begin{loglogaxis}[
			xmode=linear,
			ymode=log,
			title = {},
			xlabel={\LARGE{$N$}},
			ylabel={\LARGE{$\displaystyle{\max_{\mu \in \Xi_{\rm{test}}}} E_{\mu}$}},
			line width=1pt,
			mark size=4pt,
			grid=both,
			minor grid style={gray!25},
			major grid style={gray!25},
			width=0.75\linewidth,
			xmin=10, 
			ymin=0.0001,
			axis background/.style={fill=white}
			]
			]
			\addplot[ % <-- plot options
			%only marks,
			dashed,
			mark=o,
			mark options={solid,color = black},
			]  %
			table{dat/par_pb/errtest.dat};
			\addlegendentry{\large{POD-Greedy}}
			\addplot[ % <-- plot options
			%only marks,
			mark=triangle,
			mark options={solid,color = black},
			color=black,
			]  %
			table{dat/par_pb/errtest_stronggreedy.dat};
			\addlegendentry{\large{strong POD-Greedy}};
		\end{loglogaxis}
	\end{tikzpicture}
}
\subfloat[HA-POD Greedy]{
		\label{fig:HAPOD-predictions}
	\begin{tikzpicture}[scale=.6]
		\begin{loglogaxis}[
			xmode=linear,
			ymode=log,
			title = {},
			xlabel={\LARGE{$N$}},
			ylabel={\LARGE{$\displaystyle{\max_{\mu \in \Xi_{\rm{test}}}} E_{\mu}$}},
			line width=1pt,
			mark size=4pt,
			grid=both,
			minor grid style={gray!25},
			major grid style={gray!25},
			width=0.75\linewidth,
			xmin=10, 
			ymin=0.0001,
			%axis background/.style={fill=white}
			]
			]
			\addplot[ % <-- plot options
			%only marks,
			dashed,
			mark=o,
			mark options={solid,color = black},
			]  %
			table{dat/par_pb/hapod/errtest.dat};
			\addlegendentry{\large{POD-Greedy}}
			\addplot[ % <-- plot options
			%only marks,
			mark=triangle,
			mark options={solid,color = black},
			color=black,
			]  %
			table{dat/par_pb/hapod/errtest_stronggreedy.dat};
			\addlegendentry{\large{strong POD-Greedy}};
		\end{loglogaxis}
	\end{tikzpicture}
}
\caption{Out-of-sample performance of the ROM parametric problem obtained using the POD-Greedy algorithm. Comparison with strong POD Greedy.}
\label{fig:predictions}
\end{figure}

In Figure \ref{fig:predictions}, we assess out-of-sample performance of the proposed method. More precisely, we show the behaviour of the maximum relative error \eqref{eq:relative_L2error} over the test set $\displaystyle \max_{\mu \in \Xi_{\rm test}}E_{\mu}$ for both H-POD Greedy and HA-POD Greedy. To provide a relevant benchmark, we compare results with the H-POD Greedy and HA-POD Greedy algorithms based on the exact errors (strong POD-Greedy). For this particular example, we observe that the proposed method is effective to generate accurate ROMs: in particular, the Greedy procedures based on the time-averaged error indicator are comparable in terms of performance with the corresponding 
strong POD-Greedy algorithms.
 
\section{Conclusions}
\label{sec:conclusions}
In this work, we developed and numerically validated a model order reduction procedure for a class of problems in nonlinear mechanics, and we successfully applied it to a two-dimensional parametric THM problems that arises in radio-active waste management. We proposed a time-averaged error indicator to drive the offline Greedy sampling, and an empirical quadrature procedure to reduce offline costs.

We aim to extend the approach in several directions. First, we wish to apply our method to other problems of the form \eqref{eq:cont_pb}, to demonstrate the generality of the approach and its relevance for continuum mechanics applications. Second, we wish to combine our approach with domain decomposition methods (\cite{bergmann2018zonal,kaulmann2011new,huynh2013static}) to deal with more complex parametrizations and topological changes. Towards this end, we wish to devise effective localised training methods to reduce offline costs and domain decomposition strategies to glue together the solution in different components of the domain.

\section*{Acknowledgements}
The authors acknowledge the financial support of Andra (National Agency for Radioactive Waste Management) and thank Dr. Marc Leconte and Dr. Antoine Pasteau (Andra) for fruitful discussions.

 \newpage
\bibliographystyle{abbrv}   
\bibliography{all_refs}

\begin{thebibliography}{10}

\bibitem{barrault2004empirical}
M.~Barrault, Y.~Maday, N.~C. Nguyen, and A.~T. Patera.
\newblock An empirical interpolation method: application to efficient
  reduced-basis discretization of partial differential equations.
\newblock {\em Comptes Rendus Mathematique}, 339(9):667--672, 2004.

\bibitem{bergmann2009enablers}
M.~Bergmann, C.-H. Bruneau, and A.~Iollo.
\newblock Enablers for robust {POD} models.
\newblock {\em Journal of Computational Physics}, 228(2):516--538, 2009.

\bibitem{bergmann2018zonal}
M.~Bergmann, A.~Ferrero, A.~Iollo, E.~Lombardi, A.~Scardigli, and H.~Telib.
\newblock A zonal {G}alerkin-free {POD} model for incompressible flows.
\newblock {\em Journal of Computational Physics}, 352:301--325, 2018.

\bibitem{berkooz1993proper}
G.~Berkooz, P.~Holmes, and J.~Lumley.
\newblock The proper orthogonal decomposition in the analysis of turbulent
  flows.
\newblock {\em Annual review of fluid mechanics}, 25(1):539--575, 1993.

\bibitem{brand2003fast}
M.~Brand.
\newblock Fast online {SVD} revisions for lightweight recommender systems.
\newblock In {\em Proceedings of the 2003 SIAM international conference on data
  mining}, pages 37--46. SIAM, 2003.

\bibitem{carlberg2013gnat}
K.~Carlberg, C.~Farhat, J.~Cortial, and D.~Amsallem.
\newblock The {GNAT} method for nonlinear model reduction: effective
  implementation and application to computational fluid dynamics and turbulent
  flows.
\newblock {\em Journal of Computational Physics}, 242:623--647, 2013.

\bibitem{chapman2017accelerated}
T.~Chapman, P.~Avery, P.~Collins, and C.~Farhat.
\newblock Accelerated mesh sampling for the hyper reduction of nonlinear
  computational models.
\newblock {\em International Journal for Numerical Methods in Engineering},
  109(12):1623--1654, 2017.

\bibitem{chaturantabut2010nonlinear}
S.~Chaturantabut and D.~C. Sorensen.
\newblock Nonlinear model reduction via discrete empirical interpolation.
\newblock {\em SIAM Journal on Scientific Computing}, 32(5):2737--2764, 2010.

\bibitem{farhat2015structure}
C.~Farhat, T.~Chapman, and P.~Avery.
\newblock Structure-preserving, stability, and accuracy properties of the
  energy-conserving sampling and weighting method for the hyper reduction of
  nonlinear finite element dynamic models.
\newblock {\em International Journal for Numerical Methods in Engineering},
  102(5):1077--1110, 2015.

\bibitem{fick2018stabilized}
L.~Fick, Y.~Maday, A.~T. Patera, and T.~Taddei.
\newblock A stabilized pod model for turbulent flows over a range of reynolds
  numbers: Optimal parameter sampling and constrained projection.
\newblock {\em Journal of Computational Physics}, 371:214 -- 243, 2018.

\bibitem{edfreport2009}
S.~Granet.
\newblock Mod{\'e}lisations {THHM}. {G}{\'e}néralit{\'e}s et algorithmes.
\newblock {https://www.code-aster.org/V2/doc/v9/fr/man{\_}r/r7/r7.01.10.pdf},
  2009.

\bibitem{haasdonk2013convergence}
B.~Haasdonk.
\newblock Convergence rates of the {POD}-greedy method.
\newblock {\em ESAIM: Mathematical Modelling and Numerical Analysis},
  47(3):859--873, 2013.

\bibitem{haasdonk2017reduced}
B.~Haasdonk.
\newblock Reduced basis methods for parametrized pdes--a tutorial introduction
  for stationary and instationary problems.
\newblock {\em Model reduction and approximation: theory and algorithms},
  15:65, 2017.

\bibitem{haasdonk2008reduced}
B.~Haasdonk and M.~Ohlberger.
\newblock Reduced basis method for finite volume approximations of parametrized
  linear evolution equations.
\newblock {\em ESAIM: Mathematical Modelling and Numerical Analysis},
  42(2):277--302, 2008.

\bibitem{hesthaven2016certified}
J.~S. Hesthaven, G.~Rozza, and B.~Stamm.
\newblock {\em Certified reduced basis methods for parametrized partial
  differential equations}.
\newblock Springer, 2016.

\bibitem{himpe2018hierarchical}
C.~Himpe, T.~Leibner, and S.~Rave.
\newblock Hierarchical approximate proper orthogonal decomposition.
\newblock {\em SIAM Journal on Scientific Computing}, 40(5):A3267--A3292, 2018.

\bibitem{huynh2013static}
D.~B.~P. Huynh, D.~J. Knezevic, and A.~T. Patera.
\newblock A static condensation reduced basis element method: approximation and
  a posteriori error estimation.
\newblock {\em ESAIM: Mathematical Modelling and Numerical Analysis},
  47(1):213--251, 2013.

\bibitem{kaulmann2011new}
S.~Kaulmann, M.~Ohlberger, and B.~Haasdonk.
\newblock A new local reduced basis discontinuous {G}alerkin approach for
  heterogeneous multiscale problems.
\newblock {\em Comptes Rendus Mathematique}, 349(23-24):1233--1238, 2011.

\bibitem{larion2020building}
Y.~Larion, S.~Zlotnik, T.~J. Massart, and P.~D{\'\i}ez.
\newblock Building a certified reduced basis for coupled
  thermo-hydro-mechanical systems with goal-oriented error estimation.
\newblock {\em Computational mechanics}, 66(3):559--573, 2020.

\bibitem{lawson1974solving}
C.~L. Lawson and R.~J. Hanson.
\newblock {\em Solving least squares problems}, volume 161.
\newblock Siam, 1974.

\bibitem{leuschner2017reduced}
M.~Leuschner and F.~Fritzen.
\newblock Reduced order homogenization for viscoplastic composite materials
  including dissipative imperfect interfaces.
\newblock {\em Mechanics of Materials}, 104:121--138, 2017.

\bibitem{miled2013priori}
B.~Miled, D.~Ryckelynck, and S.~Cantournet.
\newblock A priori hyper-reduction method for coupled
  viscoelastic--viscoplastic composites.
\newblock {\em Computers \& Structures}, 119:95--103, 2013.

\bibitem{quarteroni2015reduced}
A.~Quarteroni, A.~Manzoni, and F.~Negri.
\newblock {\em Reduced basis methods for partial differential equations: an
  introduction}, volume~92.
\newblock Springer, 2015.

\bibitem{riffaud2021dgdd}
S.~Riffaud, M.~Bergmann, C.~Farhat, S.~Grimberg, and A.~Iollo.
\newblock The {DGDD} method for reduced-order modeling of conservation laws.
\newblock {\em Journal of Computational Physics}, 437:110336, 2021.

\bibitem{rozza2007reduced}
G.~Rozza, D.~B.~P. Huynh, and A.~T. Patera.
\newblock Reduced basis approximation and a posteriori error estimation for
  affinely parametrized elliptic coercive partial differential equations.
\newblock {\em Archives of Computational Methods in Engineering},
  15(3):229–275, 2007.

\bibitem{ryckelynck2009hyper}
D.~Ryckelynck.
\newblock Hyper-reduction of mechanical models involving internal variables.
\newblock {\em International Journal for Numerical Methods in Engineering},
  77(1):75--89, 2009.

\bibitem{sirovich1987turbulence}
L.~Sirovich.
\newblock Turbulence and the dynamics of coherent structures. {I.} {C}oherent
  structures.
\newblock {\em Quarterly of applied mathematics}, 45(3):561--571, 1987.

\bibitem{taddei2019offline}
T.~Taddei.
\newblock An offline/online procedure for dual norm calculations of
  parameterized functionals: empirical quadrature and empirical test spaces.
\newblock {\em Advances in Computational Mathematics}, 45(5-6):2429--2462,
  2019.

\bibitem{taddei2021discretize}
T.~Taddei and L.~Zhang.
\newblock A discretize-then-map approach for the treatment of parameterized
  geometries in model order reduction.
\newblock {\em Computer Methods in Applied Mechanics and Engineering},
  384:113956, 2021.

\bibitem{volkwein2011model}
S.~Volkwein.
\newblock Model reduction using proper orthogonal decomposition.
\newblock {\em Lecture Notes, Institute of Mathematics and Scientific
  Computing, University of Graz. see
  math.uni-konstanz.de/numerik/personen/volkwein/teaching/POD-Vorlesung.pdf},
  1025, 2011.

\bibitem{willcox2006unsteady}
K.~Willcox.
\newblock Unsteady flow sensing and estimation via the gappy proper orthogonal
  decomposition.
\newblock {\em Computers \& fluids}, 35(2):208--226, 2006.

\bibitem{yano2019discontinuous}
M.~Yano.
\newblock {Discontinuous Galerkin reduced basis empirical quadrature procedure
  for model reduction of parametrized nonlinear conservation laws}.
\newblock {\em Advances in Computational Mathematics}, pages 1--34, 2019.

\bibitem{yano2019lp}
M.~Yano and A.~T. Patera.
\newblock An {LP} empirical quadrature procedure for reduced basis treatment of
  parametrized nonlinear {PDE}s.
\newblock {\em Computer Methods in Applied Mechanics and Engineering},
  344:1104--1123, 2019.

\end{thebibliography}
 
\end{document}